\documentclass[11pt,twoside,a4paper]{article}
\usepackage{}
\usepackage{mathrsfs}
\usepackage{amsmath}
\usepackage{amsfonts}
\usepackage{amssymb}
\usepackage{amsthm}
\usepackage{graphicx,subfigure,color}
\usepackage{amsmath,amssymb}
\usepackage{indentfirst,latexsym,bm}
\usepackage{tipa}
\usepackage[mathscr]{eucal}
\begin{document}

\title{Exponential Stability Estimate of Symplectic Integrators for Integrable Hamiltonian Systems}

\author{Zhaodong Ding\footnote{School of Mathematical Sciences, Inner Mongolia University, Hohhot 010021, China;
e-mail: dingzhaodong@amss.ac.cn.} \qquad Zaijiu
Shang\footnote{Institute of Mathematics, Academy of Mathematics and
System Sciences, Chinese Academy of Sciences, Beijing 100080, P. R.
China; e-mail: zaijiu@amss.ac.cn.} \qquad Bo Xie\footnote{South China University of Technology,
Guangzhou, P. R. China; e-mail: xieb@scut.edu.cn.}}

\date{}
\maketitle

\begin{abstract}
 We prove a Nekhoroshev-type theorem for nearly integrable symplectic map. As an application of the theorem, we obtain the exponential stability symplectic algorithms. Meanwhile, we can get the bounds for the perturbation, the variation of the action variables, and the exponential time respectively. These results provide a new insight into the nonlinear stability analysis of symplectic algorithms. Combined with our previous results on the numerical KAM theorem for symplectic algorithms (2018), we give a more complete characterization on the complex nonlinear dynamical behavior of symplectic algorithms.
\end{abstract}



\section{Introduction}

\setcounter{equation}{0}

After the pioneering work of Channel (1983), Feng Kang (1985, 1986) and Ruth (1983), the symplectic integrator has become a widely interested subject on the problem of numerically solving Hamiltonian systems. Extensive computer experimentation, by some typical models of Hamiltonian systems,
has shown the overwhelming superiority of symplectic algorithms over the
conventional non-symplectic ones, especially in simulating the global and
structural dynamic behavior of the systems (e.g. see \cite{Haire} and \cite{FQ}). The symplectic algorithm, which applied to integrable Hamiltonian system, may be characterized as a perturbation of the phase flow of the integrable system. Here the smallness of the perturbation is described by the
time-step size of the algorithm which also enters into the frequency map of the integrable system. Therefore numerical stability problem arises.

We consider a nearly integrable symplectic map $\mathcal
{C}:(I,\theta)\rightarrow (\hat{I},\hat{\theta})$ generated by an
analytical function in action-angle variables of the form
$H(\hat{I},\theta)=H_0(\hat{I})+h(\hat{I},\theta)$ which is defined
on $\mathcal {G}\times T^n$, where actions space $\mathcal {G}$ is
an open and bounded domain of $\mathbb{R}^n$ and $h$ is a small
perturbation, of size $\epsilon$. The symplectic map $\mathcal {C}$
is given implicitly by
\begin{equation*} \label{eq:1}
 \begin{aligned}
         \hat{I} &= I-\partial_2h(\hat{I},\theta)~, \\
         \hat{\theta} &=\theta+\partial H_0(\hat{I})+\partial_1h(\hat{I},\theta)~.
                          \end{aligned}
                          \end{equation*}

It is clear that the symplectic map is integrable when the perturbation $h$ vanishes. In
this case the dynamics of the symplectic map is trivial. The action variables remain constant
for all iterative times and the angle variables vary linearly with respect to the iteration steps.
However, a perturbed symplectic map may generate very complicated dynamics in general
when the perturbation $h$ is non-vanishing. It is believed that there exists unstable motions.
For example, Arnold diffusion may take place if the degrees of freedom are higher than two
in mapping case \cite{Arnold-2}.

By KAM theorem, nevertheless, most of motions are perpetually stable for sufficiently
smooth nearly integrable Hamiltonian systems, and these stable ones form a Cantor set
with large measure in the phase space \cite{Arnold,Shang2000}. On the other hand, all motions
are exponentially stable if the systems are analytic and the unperturbed integrable part of
the Hamiltonian satisfies the so-called steepness condition by Nekhoroshev theorem. The result was first proved by Nekhoroshev for nearly integrable Hamiltonian systems and
are also called effective stability in literature. More precisely, he
proved for sufficiently small $\epsilon$ and for all initial values
$I(0)$ in action space, one has
$$|I(t)-I(0)|_2\leq C\epsilon^b\qquad \mbox{for}\quad |t|\leq T_0\exp(\epsilon^{-a})~,$$
with constants $C,~T_0,~a$, and $b$, provided the steepness conditions
are fulfilled for $H_0$. Further improvements on the stability exponents $a$ and $b$
were made by P\"{o}schel \cite{Poschel}, Lochak \cite{Lochak-1,Lochak-2} and Bounemoura \cite{AB}.

In his paper, Nekhoroshev conjectured similar result will hold for the nearly integrable symplectic map. Later, Kuksin and P\"{o}schel
\cite{Kuksin-Poschel} gave a proof about the exponential stability
of nearly integrable symplectic maps by proving the existence of a
non-autonomous analytic Hamiltonian system interpolating a
symplectic map and applying the Nekhoroshev theorem of the
Hamiltonian case. However, the proof is an existence one and there
is not explicit estimate about the small perturbation. In 2004,
Guzzo \cite{Guzzo} proposed a direct proof of the Nekhoroshev
theorem for nearly integrable symplectic maps. This result is very important and valuable, but his estimate
for the perturbation is not suitable to be applied to the small
twist problem and thus the symplectic integrator. Here we get a different estimate of the perturbation so
that it can be applied to the small twist maps and symplectic
integrators. In \cite{XieBo}, Nekhoroshev stability of symplectic
algorithms has been considered. But there the Nekhoroshev stability can be
obtained only when the order of the algorithm is greater than 2. In this article we completely solve this problem.

In the present paper, we obtain an exponential stability result by
construction, and provide explicit estimates of the involved
quantities for the nearly integrable symplectic map, then apply it
to the cases of small twist maps and symplectic integrators.
Following the original idea of Nekhoroshev \cite{Nekhoroshev}, the
proof of the main theorem is divided into three steps. First, normal forms of
the nearly integrable symplectic map are constructed on some subdomains of phase space that are
known as resonant blocks. Second, these normal forms
lead to stability estimates on the corresponding subdomains. The third step is to give a
geometric construction which guarantees the whole action space can be completely covered
by such resonant blocks.

\subsection{Notations}

We introduce some notations used in this paper, most of which are
from \cite{Delshams-Gutierrez}. Given $\rho=(\rho_1,\rho_2)\geq 0$
(i.e. $\rho_j\geq 0$, $j=1,2$), we first introduce the sets
$$\mathcal {V}_{\rho_1}(G):=\{I\in \mathbb{C}^n~|~|I-I'|_2\leq \rho_1~\mbox{for some}
~I'\in G\},$$ and
$$\mathcal {W}_{\rho_2}(T^n):=\{\theta\in \mathbb{C}^n/(2\pi\mathbb{Z}^n)~|~\texttt{Re}\theta\in T^n,~|\texttt{Im}\,\theta|_{\infty}\leq \rho_2\},$$
where $|\cdot|_2$ and $|\cdot|_{\infty}$ denote, respectively, the
Euclidean norm and the maximum norm for vectors; $\texttt{Re}\theta$ and $\texttt{Im}\,\theta$ denote the real part and the imaginary part of $\theta$ respectively. Then define
$$\mathscr{D}_{\rho}(G):=\mathcal {V}_{\rho_1}(G)\times \mathcal {W}_{\rho_2}(T^n)\,.$$

Several kinds of norms are used along this paper. First, we consider
functions of the $n$ action variables. Given a (real or complex)
function $f(I)$, defined on a complex neighborhood $\mathcal
{V}_{\eta}(G)$, we introduce the supremum norm
$$|f|_{G,\,\eta}:=\sup_{I\in \mathcal
{V}_{\eta}(G)}|f(I)|\,,\qquad |f|_{G}:=|f|_{G,\,0}\,.$$ In this way,
the subscript $\eta$ is removed from the notation if $\eta=0$. This
remark applies throughout this section.

In an analogous way, we consider the supremum norm for vector-valued
functions, i.e., vectorfields. Given $F: \mathcal
{V}_{\eta}(G)\rightarrow\mathbb{C}^n$ and $1\leq p\leq \infty$, we
define
$$|F|_{G,\,\eta,\,p}:=\sup_{I\in \mathcal
{V}_{\eta}(G)}|F(I)|_p\,,\qquad
|F|_{G,\,\eta}:=|F|_{G,\,\eta,\,2}\,.$$ In this definition,
$|\cdot|_p$ means the p-norm for vectors in $\mathbb{C}^n$, i.e.
$|v|_p=(\sum_{j=1}^n|v_j|^p)^{1/p}$ for $1\leq p<\infty$, and
$|v|_{\infty}=\max_{1\leq j\leq n}|v_j|$\,.

Next we consider functions of the action-angle variables. For a
given complex function $f(I,\theta)$ ($2\pi$-periodic in $\theta$)
defined on the neighborhood
$\mathscr{D}_{\rho}(G),~\rho=(\rho_1,\,\rho_2)\geq 0$, we may
consider its supremum norm
$$|f|_{G,\,\rho}:=\sup_{(I,\,\theta)\in \mathscr{D}_{\rho}(G)}|f(I, \theta)|\,.$$
But if $f$ is analytic on (a neighborhood of ) the set
$\mathscr{D}_{\rho}(G)$, we may define an exponentially weighted
norm in terms of the Fourier series of $f$. Writing
$f(I,\theta)=\sum_{k\in \mathbb{Z}^n}f_k(I)e^{ik\cdot\theta}$, we
introduce
$$||f||_{G,\,\rho}:=\sum_{k\in \mathbb{Z}^n}|f_k|_{G,\,\rho_1}\cdot e^{|k|_1\rho_2}\,. $$
Note that $|f|_{G,\,\rho}\leq ||f||_{G,\,\rho}$\,.

Exactly in the same way as before we may extend the definitions of
the norms to the case of vector-valued functions. Given
$F:\mathscr{D}_{\rho}(G)\rightarrow\mathbb{C}^n$ and $1\leq p\leq
\infty$, and writing $F(I,\theta)=\sum_{k\in
\mathbb{Z}^n}F_k(I)e^{ik\cdot\theta}$, where $F_k:\mathcal
{V}_{\rho_1}(G)\rightarrow\mathbb{C}^n$, we define
$$||F||_{G,\,\rho,\,p}:=\sum_{k\in \mathbb{Z}^n}|F_k|_{G,\,\rho_1,\,p}\cdot e^{|k|_1\rho_2}\,, \qquad ||F||_{G,\,\rho}:=||F||_{G,\,\rho,\,2}\,.$$

The Cauchy estimates about the Fourier norms are provided by
\cite{Poschel}. That is, if $f$ analytic on $\mathscr{D}_{\rho}(G)$,
for $0<\delta<\rho$ one has
$$\Big|\Big|\frac{\partial f}{\partial I}\Big|\Big|_{G,\,(\rho_1-\delta_1,\,\rho_2),\,\infty}\leq \frac{1}{\delta_1}||f||_{G,\,\rho}\,,\qquad
\Big|\Big|\frac{\partial f}{\partial
\theta}\Big|\Big|_{G,\,(\rho_1,\,\rho_2-\delta_2),\,1}\leq
\frac{1}{e\delta_2}||f||_{G,\,\rho}\,.$$

Finally, for $Df=(\partial f/\partial I,\,\partial f/\partial
\theta)$ we introduce the vectorfield norm
$$||Df||_{G,\,\rho,\,c}:=\max\Big(\Big|\Big|\frac{\partial f}{\partial
\theta}\Big|\Big|_{G,\,\rho,\,1},\,c\Big|\Big|\frac{\partial
f}{\partial I}\Big|\Big|_{G,\,\rho,\,\infty}\Big),$$ where $c>0$ is
a parameter to be fixed in subsequent sections.

\subsection{Main Result}

The main theorem for nearly integrable symplectic map can be stated
as follows.\\[-5pt]

\textbf{Theorem 1.} Let
$H(\hat{I},\theta)=H_0(\hat{I})+h(\hat{I},\theta)$ be analytic in
$\mathscr{D}_{\sigma}(\mathcal {G})$, where $\mathcal {G}$ is an
open bounded domain of $\mathbb{R}^n$, and
$\sigma=(\sigma_1,\sigma_2)$ is positive. Let
$\omega(\hat{I})=\partial H_0(\hat{I})$ satisfying
\begin{equation}\label{eq:30}
m|I_1-I_2|_2\leq |\omega(I_1)-\omega(I_2)|_2\leq M|I_1-I_2|_2
\end{equation}
for $I_1,I_2\in \mathcal {V}_{\sigma_1}(\mathcal {G})$ with positive
constants $m, \,M$.

Assume
$$\epsilon=||h||_{\mathcal {G},\,\sigma}\leq \min(\epsilon_0,\,\sigma_2^b)\,,$$
where
$$\epsilon_0=\frac{M^2\cdot\sigma_1^4}{\pi^2(21n+30)^2\big[(\frac{2M}{m})^n\cdot n!\big]^4}\,,$$
$b=2(n^2+n+2)$ and $\sigma_1< \frac{1}{4M}$\,.

Then the symplectic map $\mathcal {C}$ which generated by
$H(I',\theta)$ satisfies
$$|I(t)-I_0|_2\leq \triangle \quad \mbox{for}\quad |t|\leq T \quad \mbox{and}\quad I_0\in \mathcal {G}-\triangle\,,$$
where $(I(t),\theta(t))=\mathcal {C}^t(I_0,\theta_0)$ ($t$ is viewed
as iterative times),
\begin{equation*} \label{eq:1}
 \begin{aligned}
         \triangle &= c_0\epsilon^{\frac{1}{b}}\quad&\mbox{with}\quad &c_0=\frac{8nM}{3m}(3n+2)\sigma_1\,,\\
         T &=T_0\epsilon^{-\frac{3}{4}}e^{c_1\epsilon^{-\frac{1}{b}}}\quad&\mbox{with}\quad
         &T_0=\frac{(\frac{m}{2M})^n\sigma_1\sigma_2}{2^7n!}\quad\mbox{and}\quad
         c_1=\frac{\sigma_2}{24}\,.
                          \end{aligned}
                          \end{equation*}\\[-5pt]

The rest of this paper is devoted to the proof of the main theorem.
In Section 2, the analytic part is presented which concerns the
construction of the normal form with exponentially small remainder
on resonant blocks. The geometric part, in Section 3, concerns the
covering of the whole action space $\mathcal {G}$ by a family of
resonant blocks. In Section 4, the proof of
Theorem 1 is finished through making choices for the free
parameters. Finally, as applications of the main theorem, we consider the stability of the small twist
type symplectic map and symplectic integrator which applied to the
integrable Hamilton system in
Section 5.

\section{The Analytic Part}

At first, we transform the mapping $\mathcal {C}$ by the partial
coordinates stretching $\mathcal {W}_{\gamma}:(x,y)\rightarrow
(I,\theta)=(\gamma x,y)$\,, and obtain a new mapping
$T_{\gamma}=\mathcal {W}_{\gamma}^{-1}\circ \mathcal {C}\circ
\mathcal {W}_{\gamma} :(x,y)\rightarrow (\hat{x},\hat{y})$ to be
defined in the new phase space $\mathcal {G}_{\gamma}\times T^n$ by

\begin{equation*} \label{eq:2}
\left\{ \begin{aligned}
         \hat{x} &= x-\partial _2 F(\hat{x},y) \\
                  \hat{y}&=y+\partial _1 F(\hat{x},y)
                          \end{aligned} \right.
                          \end{equation*}
where
$$F(x,y)=F_0(x)+f(x,y)$$
is well defined on $\mathcal {G}_{\gamma}\times T^n$ with
$$F_0=\gamma^{-1}H_0(\gamma x),\quad f(x,y)=\gamma^{-1}h(\gamma x,y)$$
and
$$\mathcal {G}_{\gamma}=\{x\in \mathbb{R}^n~|~\gamma x\in\mathcal
{G}\}\,.$$

For the time being, $\gamma$ is considered as a free parameter.
$F(x,y)$ is real analytic in $\mathcal
{D}_{\widetilde{\sigma}}(\mathcal {G}_{\gamma})$, where
$\widetilde{\sigma}=(\widetilde{\sigma}_1,\widetilde{\sigma}_2)$
with $\widetilde{\sigma}_1=\gamma^{-1}\sigma_1$ and
$\widetilde{\sigma}_2=\sigma_2$. Accordingly, the frequency map of
the integrable mapping associated to the generating function $F_0$
turns into $\widetilde{\omega}(x)=\partial F_0(x)$ and the condition
satisfied by the map $\widetilde{\omega}$ turns out to be
\begin{equation}
\gamma m|x_1-x_2|_2\leq
|\widetilde{\omega}(x_1)-\widetilde{\omega}(x_2)|_2 \leq \gamma
M|x_1-x_2|_2
\end{equation}
for $x_1,x_2\in \mathcal {V}_{\widetilde{\sigma}_1}(\mathcal
{G}_{\gamma})$\,. In addition, we have
$$\tilde{\epsilon}:=||f||_{\mathcal {G}_{\gamma},\,\widetilde{\sigma}}=\gamma^{-1}||h||_{\mathcal {G},\,\sigma}=\gamma^{-1}\epsilon\,.$$

From now on, we fix $\gamma=\beta M^{-1}$, where $\beta>0$
is a parameter to be determined later. Denoting $\mu=\frac{m}{M}$\,,
we have
\begin{equation}\label{eq:3}
\beta\mu |x_1-x_2|_2\leq
|\widetilde{\omega}(x_1)-\widetilde{\omega}(x_2)|_2 \leq \beta
|x_1-x_2|_2
\end{equation}

In order to transform the generating function $F(x,y)=F_0(x)+f(x,y)$
of $T_{\gamma}$ to the normal form, we need seek for
a suitable canonical transformation $\Phi$, that is constructed
iteratively as a product of the successive approximately identical canonical transformation
$\Phi^{(1)},~\Phi^{(2)},\ldots$. In doing that, one meets the small
denominators $1-e^{ik\cdot\widetilde{\omega}(x)}$, which in general
vanish in a dense subset of $\mathcal {G}_{\gamma}$\,. These
resonances are given by the equations
$k\cdot\widetilde{\omega}(x)+2\pi l=0$, for $k\in \mathbb{Z}^n,~l\in
\mathbb{Z}$. As usual, the small denominators
$1-e^{ik\cdot\widetilde{\omega}(x)}$ should been excluded in the
process of constructing the normal form. However, it is not necessary to take care
of all the small denominators. We can do that in a certain subdomain which can guarantee
$|1-e^{ik\cdot\widetilde{\omega}(x)}|\geq \alpha$ for $|k|\leq K$ and $k\notin\mathcal {M}$,
where $\alpha$ and $K$ are independent parameters and $\mathcal {M}$ is a given sublattice of $\mathbb{Z}^n$.

Next we give the related notations in detail. Let $\mathcal {M}$ be a sublattice of $\mathbb{Z}^n$. We only
consider the maximal ones, that is not properly contained in any
other sublattice of the same dimension. A maximal sublattice
$\mathcal {M}$ with dim$\mathcal {M}=r$, is said $K$-lattice if it
admits a basis $(k^{(1)},\ldots,k^{(r)})$ satisfying
$|k^{(j)}|_1\leq K$ for $j\leq r$ holds, and such a basis will be called a $K$-basis (\cite{Benettin}). A function $g(x,y)$ is said
to be in a normal form with respect to $\mathcal {M}$ of degree $K$
if its Fourier series expansion in the angular variables is
restricted to the form $g(x,y)=\sum_{k\in \mathcal {M},\,|k|_1\leq
K}g_k(x)e^{ik\cdot y}$. We express this by writing $g\in \mathcal
{R}(\mathcal {M},K)$. Note that a function is in the normal form
with respect to the trivial modulo $\mathcal {M}=0$ if it does not
depend on the angular variables.

We restrict ourselves to a subset $G\subset \mathcal {G}_{\gamma}$,
where the frequency vectors $\widetilde{\omega}(x)$ is allowed to
satisfy some resonance relations corresponding to a fixed sublattice
$\mathcal {M}$, but a neighborhood of all other resonances of order
less than or equal to $K$ are excluded. In precise, a subset
$G\subset \mathcal {G}_{\gamma}$ is said to be
$\alpha,K$-nonresonant modulo $\mathcal {M}$ if
$$|1-e^{ik\cdot\widetilde{\omega}(x)}|\geq \alpha\quad \mbox{for all}~\, k\in \mathbb{Z}^n_K\backslash \mathcal {M}~\, \mbox{and}~\, x\in G\,,$$
where $\mathbb{Z}^n_K:=\{k\in\mathbb{Z}^n~|~|k|_1\leq K\}$.

A special situation arises when $\mathcal {M}$ is the trivial
sublattice $\Theta$ of $\mathbb{Z}^n$ containing only 0. In this
case the set $G$ is said to be completely $\alpha,K$-nonresonant. In
the corresponding normal form, $g$ is independent of the angle
variables. Naturally, the analysis of this case is simpler than in
the presence of resonances.

The nonresonance condition on the set $G$ can be extended to a
complex neighborhood of small enough radius $\rho_1$.\\[-5pt]

\textbf{Lemma 2.1} Let $F_0(x)$ be a real analytic function in
$\mathcal {V}_{\rho_1}(G)$, and let $\widetilde{\omega}=\partial
F_0$. Assume that $G$ is $\alpha,K$-nonresonant modulo $\mathcal
{M}$ $(0<\alpha\leq 1)$, and $\widetilde{\omega}$ satisfies
(\ref{eq:3}). If
$$\rho_1\leq \frac{\alpha}{4K\beta}\,,$$
then $\mathcal {V}_{\rho_1}(G)$ is $\alpha/2,K$-nonresonant modulo
$\mathcal {M}$.\\[-5pt]

\textbf{Proof.} $\forall x\in\mathcal {V}_{\rho_1}(G)$, there exists
$\tilde{x}\in G$ such that $|x-\tilde{x}|_2\leq \rho_1$ by
definition. Because
\begin{equation*}
\begin{split}
\big|1-e^{ik\cdot\widetilde{\omega}(x)}\big|&\geq  \big|1-e^{ik\cdot\widetilde{\omega}(\tilde{x})}\big|-\big|e^{ik\cdot\widetilde{\omega}(x)}-e^{ik\cdot\widetilde{\omega}(\tilde{x})}\big| \\
 &\geq \alpha-\big|e^{ik\cdot\widetilde{\omega}(x_{*})}<k,\,\widetilde{\omega}(\tilde{x})-\widetilde{\omega}(x)>\big|\,,
 \end{split}
 \end{equation*}
where for the last inequality we have used the mean value theorem
with $x_{*}\in \mathcal {V}_{\rho_1}(G)$. Thus, there exists $x'\in
G$, s.t. $|x_{*}-x'|_2\leq \rho_1$ and $|\bar{x}_{*}-\bar{x}'|_2\leq
\rho_1$, where $\bar{x}$ is the complex conjugation of $x$. Noting
that $x'=\bar{x}'$, we have $|x_{*}-\bar{x}_{*}|_2\leq 2\rho_1$.
Therefore,
$$\big|e^{ik\cdot\widetilde{\omega}(x_{*})}\big|\leq e^{|k|_1\cdot|\texttt{Im}\,\widetilde{\omega}(x_{*})|_2}\leq
e^{K\cdot\frac{|\widetilde{\omega}(x_{*})-\widetilde{\omega}(\bar{x}_{*})|_2}{2}}\leq
e^{K\beta\rho_1}\leq 2$$ and
$\big|e^{ik\cdot\widetilde{\omega}(x_{*})}<k,\,\widetilde{\omega}(\tilde{x})-\widetilde{\omega}(x)>\big|\leq
2K\beta\rho_1\leq \frac{\alpha}{2}$\,. That completes our proof.
\qed

Next we state the Normal Form Lemma as follows.\\[-5pt]

\textbf{Lemma 2.2 (Normal Form Lemma)} Let $\mathcal
{M}\subseteq\mathbb{Z}^n$ be a $K$-lattice, $0<\alpha\leq 1$ and $\beta\leq \alpha$.
$F(x,y)=F_0(x)+f(x,y)$ is analytic in $\mathcal {D}_{\rho}(G)$,
$\rho=(\rho_1,\rho_2)$. Suppose that $\mathcal {V}_{\rho_1}(G)$ is
$\alpha,K$-nonresonant modulo $\mathcal {M}$ and the frequency map
$\widetilde{\omega}$ satisfies (\ref{eq:3}). If
\begin{equation}\label{eq:5}
||Df||_{G,\,\rho,\,c}\leq \frac{\alpha\,\rho_1}{\widetilde{C}AK\rho_2}\,,
\end{equation}
where $\widetilde{C}=21n+30,~A=1+\frac{\beta
c}{\alpha}e^{K\beta\rho_1}$ and $\rho_1\leq
\min(\frac{\alpha}{4K\beta},1)$, then there exists a real analytic
canonical transformation $\Phi :\mathcal
{D}_{\frac{\rho}{2}}(G)\rightarrow \mathcal {D}_{\rho}(G)$ such that
the conjugate symplectic map $T_{\gamma}'=\Phi^{-1}\circ
T_{\gamma}\circ\Phi:\mathcal {D}_{\frac{\rho}{2}}(G)\rightarrow
\mathcal {D}_{\rho}(G)$ is generated by the analytic function
$F'=F_0+Z+R$ with $Z\in \mathcal {R}(\mathcal {M},K)$. Moreover, the
following hold\\[-6pt]

1)\quad$||DZ||_{G,\,\frac{\rho}{2},\,c}\leq
2||Df||_{G,\,\rho,\,c}$\,.\\[-7pt]

2)\quad $||DR||_{G,\,\frac{\rho}{2},\,c}\leq
3e^{-\frac{K\rho_2}{12}}||Df||_{G,\,\rho,\,c}$\,.\\[-7pt]

3)\quad $|P_{x}\Phi-id|_{G,\,\frac{\rho}{2}}\leq
\frac{\rho_1}{2^8}$\,,\\[6pt]
where $P_x$ denotes the projection onto $x$-coordinates.\\[-5pt]

In order to prove the Normal Form Lemma, we need some technique
lemmas.

\textbf{Lemma 2.3 \cite{Delshams-Gutierrez}} Let $f$ be an analytic
function in $\mathcal {D}_{\rho}(G)$. For
$0<\delta=(\delta_1,\delta_2)<\rho$ and $c>0$ given, let us denote
$$\hat{\delta}_c:=\min(\delta_1,c\delta_2)\,.$$
Then,\\[-6pt]

(a)\quad $||Df||_{G,\,\rho-\delta,\,c}\leq
\frac{c}{\hat{\delta}_c}||f||_{G,\,\rho}\,.$\\[-7pt]

(b)\quad $||D(f^{>K})||_{G,\,\rho-\delta,\,c}\leq
e^{-K\delta_2}||Df||_{G,\,\rho,\,c}\,,$ where
$f^{>K}(x,y)=\sum\limits_{|k|_1>K}f_k(x)e^{ik\cdot y}$\,.\\[-5pt]

The following lemma is similar with Lemma 6 in \cite{Guzzo}, and the
way of the proof can be found there.\\
 \textbf{Lemma 2.4 \cite{Guzzo}} Let $T_{\gamma}$ and $\Phi$ be
 symplectic maps which generated by analytic functions $F$ and
 $\chi$ respectively. $T_{\gamma},\Phi :G\times T^n\rightarrow G\times
 T^n$ defined by
 \begin{gather}
  \Phi :(a,\varphi)\rightarrow (x,y) \notag\\
  \left\{ \begin{aligned}
         x&=a+\partial_2 \chi(a,y)\\
                   y&=\varphi-\partial_1 \chi(a,y)
                           \end{aligned} \right.
\end{gather}
\begin{gather}
  T_{\gamma} :(x,y)\rightarrow (\hat{x},\hat{y}) \notag\\
  \left\{ \begin{aligned}
         \hat{x}&=x-\partial_2 F(\hat{x},y)\\
                   \hat{y}&=y+\partial_1 F(\hat{x},y)
                           \end{aligned} \right.
\end{gather}

Then one of the generating functions of the conjugate symplectic map
$T_{\gamma}'=\Phi^{-1}\circ T_{\gamma}\circ\Phi
:(a,\varphi)\rightarrow (\hat{a},\hat{\varphi})$ is
$$\widetilde{F}(\hat{a},\varphi)=a\cdot\varphi-a\cdot
y+\hat{a}\cdot\hat{y}-\hat{a}\cdot \varphi+\hat{x}\cdot
y-\hat{x}\cdot\hat{y}+F(\hat{x},y)+\chi(\hat{a},\hat{y})-\chi(a,y)\,,$$
where variables $x,\hat{x},y,\hat{y},a,\hat{\varphi}$ are functions
of the independent variables $\hat{a},\varphi$\,.\\[-5pt]

\textbf{Lemma 2.5} Let $\chi(x,y)$ be analytic in $\mathcal
{D}_{\rho}(G)$, and given positive numbers
$\delta=(\delta_1,\delta_2)<(\rho_1,\rho_2)$. Assume
$||D\chi||_{G,\,\rho,\,c}\leq \frac{\hat{\delta}_c}{2}$, then the
symplectic transformation $\Phi$ generated by the function $\chi$ is
well defined in $\mathcal {D}_{\rho-\delta}(G)$. Furthermore, one
has
$$\Phi(\mathcal {D}_{\rho-\delta}(G))\subseteq \mathcal {D}_{\rho-\frac{\delta}{2}}(G)\quad \mbox{and}\quad
\Phi^{-1}(\mathcal {D}_{\rho-\delta}(G))\subseteq \mathcal
{D}_{\rho-\frac{\delta}{2}}(G). $$

\textbf{Proof.} Let $\mathcal {A}_{\delta_1,\,\delta_2}=\{v:\mathcal
{D}_{\rho-\delta}(G)\rightarrow
\mathbb{C}^n~|~||v||_{G,\,\rho-\delta}\leq \frac{\delta_2}{2}\}$. It
is easy to know that $\mathcal {A}_{\delta_1,\,\delta_2}$ is a
closed bounded subset of a Banach space. Now consider the map $\mathcal
{F}(v)(x,y)=\frac{\partial \chi(x,y-v(x,y))}{\partial x}$\,, which
is well defined for any $v \in\mathcal {A}_{\delta_1,\,\delta_2}$
and maps the space into itself from the fact that
\begin{equation}\label{eq:4}
||\mathcal {F}||_{G,\,\rho-\delta}\leq\frac{1}{c}||D\chi||_{G,\,\rho,\,c}\leq
\frac{\hat{\delta}_c}{2c}\leq \frac{\delta_2}{2}\,.
\end{equation}
Moreover, for any $v_1,v_2 \in\mathcal {A}_{\delta_1,\,\delta_2}$,
one has
\begin{align*}
||\mathcal {F}(v_1)(x,y)-\mathcal
{F}(v_2)(x,y)||_{G,\,\rho-\delta}&=\Big|\Big|\frac{\partial
\chi(x,y-v_1(x,y))}{\partial x}-\frac{\partial
\chi(x,y-v_2(x,y))}{\partial x}\Big|\Big|_{G,\,\rho-\delta} \\
&=\Big|\Big|\frac{\partial^2 \chi(x,y^{*})}{\partial x\partial
y}\cdot(v_1-v_2)\Big|\Big|_{G,\,\rho-\delta} \\
&\leq \frac{1}{\delta_1}\Big|\Big|\frac{\partial \chi}{\partial
y}\Big|\Big|_{G,\,\rho,1}\cdot||v_1-v_2||_{G,\,\rho-\delta}\\
&\leq \frac{1}{2}||v_1-v_2||_{G,\,\rho-\delta}\,.
\end{align*}

This shows that the map $\mathcal {F}:\mathcal
{A}_{\delta_1,\,\delta_2}\rightarrow \mathcal
{A}_{\delta_1,\,\delta_2}$ is contractive. Therefore, there exists a
unique $v^{*}\in \mathcal {A}_{\delta_1,\,\delta_2}$ such that
$\mathcal {F}(v^{*})=v^{*}$ and the symplectic transformation
$\Phi:(x,y)\rightarrow (\hat{x},\hat{y})$ can be expressed
explicitly in the form
\begin{equation*}
\left\{ \begin{aligned}
         \hat{x}&=x+\frac{\partial \chi(x,y-v^{*})}{\partial y} \\
                  \hat{y}&=y-\frac{\partial \chi(x,y-v^{*})}{\partial x}
                          \end{aligned} \right.
                          \end{equation*}
which is well defined and real analytic for any $(x,y)\in \mathcal
{D}_{\rho-\delta}(G)$. It is easy to show $\Phi(\mathcal
{D}_{\rho-\delta}(G))\subseteq \mathcal
{D}_{\rho-\frac{\delta}{2}}(G)$ by means of (\ref{eq:4}) and
$$\Big|\Big|\frac{\partial
\chi(x,y-v^{*})}{\partial
y}\Big|\Big|_{G,\,\rho-\delta,\,1}\leq
||D\chi||_{G,\,\rho,\,c}\leq \frac{\delta_1}{2}\,.$$

In similar way, we can prove that $\Phi^{-1}(\mathcal
{D}_{\rho-\delta}(G))\subseteq \mathcal
{D}_{\rho-\frac{\delta}{2}}(G)$\,. Specifically, denote the set
$\widetilde{\mathcal {A}}_{\delta_1,\,\delta_2}=\{u:\mathcal
{D}_{\rho-\delta}(G)\rightarrow
\mathbb{C}^n~|~||u||_{G,\,\rho-\delta}\leq \frac{\delta_2}{2}\}$ and
consider the map $\widetilde{\mathcal {F}}(u)(x,y)=\frac{\partial
\chi(x+u(x,y),y)}{\partial y}$\,, which is well defined in
$\widetilde{\mathcal {A}}_{\delta_1,\,\delta_2}$\,. \qed\\[-5pt]

\textbf{Lemma 2.6} Let $F(\hat{x},y)=F_0(\hat{x})+f(\hat{x},y)$ be
analytic in $\mathcal {D}_{\rho}(G)$ with
$\widetilde{\omega}(x)=\partial F_0(x)$ satisfying (\ref{eq:3}).
Given positive numbers $\delta=(\delta_1,\delta_2)<(\rho_1,\rho_2)$,
if
$$\rho_1\leq \frac{\delta_2}{6\beta}\quad \mbox{and}\quad ||Df||_{G,\,\rho,\,c}\leq \frac{\hat{\delta}_c}{2}\,,$$
then the symplectic map $T_{\gamma}$ generated by $F$ is well
defined in $\mathcal {D}_{\rho-\delta}(G)$ and
$$T_{\gamma}(\mathcal {D}_{\rho-\delta}(G))\subseteq \mathcal {D}_{\rho_1-\frac{\delta_1}{2},\,\rho_2-\frac{\delta_2}{3}}(G)\,,\quad \quad
T_{\gamma}^{-1}(\mathcal {D}_{\rho-\delta}(G))\subseteq \mathcal
{D}_{\rho_1-\frac{\delta_1}{2},\,\rho_2-\frac{\delta_2}{3}}(G). $$

\textbf{Proof.} Similar to the proof of Lemma 2.5, we know that for any
$(x,y)\in\mathcal {D}_{\rho-\delta}(G)$, there exists a unique
$\hat{x}\in \mathcal {V}_{\rho_1-\frac{\delta_1}{2}}(G)$ such that
$\hat{x}=x-\frac{\partial f(\hat{x},y)}{\partial y}$. Therefore,
$\hat{y}=y+\widetilde{\omega}(\hat{x})+\frac{\partial
f(\hat{x},y)}{\partial \hat{x}}$ is also well defined in $\mathcal
{D}_{\rho-\delta}(G)$. Moreover, when $(x,y)\in\mathcal {D}_{\rho-\delta}(G)$,
\begin{align*}
|\texttt{Im}\,(\hat{y}-y)|_{\infty}&\leq
|\texttt{Im}\,\widetilde{\omega}(\hat{x})|_{G,\,\rho_1-\delta_1/2,\,\infty}+\Big|\Big|\frac{\partial
f}{\partial x}\Big|\Big|_{G,\,\rho,\,\infty} \\
&\leq
\Big|\frac{\widetilde{\omega}(\hat{x})-\widetilde{\omega}(\bar{\hat{x}})}{2}\Big|_{G,\,\rho_1-\delta_1/2,\,\infty}+\frac{1}{c}||Df||_{G,\,\rho,\,c}\\
&\leq \frac{\beta}{2}|\hat{x}-\bar{\hat{x}}|_2+\frac{\delta_2}{2}\,.
\end{align*}
Due to $\hat{x}\in \mathcal {V}_{\rho_1-\frac{\delta_1}{2}}(G)$,
there exists $x'\in G$ such that $|\hat{x}-x'|_2 \leq \rho_1$\,.
Thus
$$|\hat{x}-\bar{\hat{x}}|_2\leq |\hat{x}-x'|_2+|\bar{x^{'}}-\bar{\hat{x}}|_2\leq 2\rho_1\leq \frac{\delta_2}{3\beta}\,.$$
Therefore, $|\texttt{Im}\,(\hat{y}-y)|_{\infty}\leq
\frac{2}{3}\delta_2\,,$ and we get
$$T_{\gamma}(\mathcal {D}_{\rho-\delta}(G))\subseteq \mathcal {D}_{\rho_1-\frac{\delta_1}{2},\,\rho_2-\frac{\delta_2}{3}}(G)\,.$$
In similar way, one can prove that $T_{\gamma}^{-1}(\mathcal
{D}_{\rho-\delta}(G))\subseteq \mathcal
{D}_{\rho_1-\frac{\delta_1}{2},\,\rho_2-\frac{\delta_2}{3}}(G)$\,.\qed\\[-5pt]

Combing Lemma 2.4, 2.5 and 2.6, we have

\textbf{Lemma 2.7} Consider the generating function
$F(\hat{x},y)=F_0(\hat{x})+f(\hat{x},y)$ and $\chi$, which are
analytic in $\mathcal {D}_{\rho}(G)$. let
$\widetilde{\omega}(x)=\partial F_0(x)$ satisfies (\ref{eq:3}).
Given positive numbers $\delta=(\delta_1,\delta_2)$ such that
$3\delta\leq \rho$\,. Suppose
$$\rho_1\leq\frac{\delta_2}{6\beta}\,,\quad||Df||_{G,\,\rho,\,c}\leq \frac{\hat{\delta}_c}{2}\quad \mbox{and}\quad ||D\chi||_{G,\,\rho,\,c}\leq \frac{\hat{\delta}_c}{2}\,.$$
Let $T_{\gamma}$ be the symplectic map generated by $F$, $\Phi$ by
$\chi$ and $T_{\gamma}'=\Phi^{-1}\circ T_{\gamma}\circ\Phi$\,. Then
$T_{\gamma}'$ and ${T_{\gamma}'}^{-1}$ are analytic and symplectic
diffemorphisms which defined in $\mathcal
{D}_{\rho_1-2\delta_1,\,\rho_2-3\delta_2}(G)$. It is
$$T_{\gamma}'(\mathcal
{D}_{\rho_1-2\delta_1,\,\rho_2-3\delta_2}(G))\subseteq \mathcal
{D}_{\rho_1-\frac{\delta_1}{2},\,\rho_2-\frac{4}{3}\delta_2}(G)\,,$$
and $${T_{\gamma}'}^{-1}(\mathcal
{D}_{\rho_1-2\delta_1,\,\rho_2-3\delta_2}(G))\subseteq \mathcal
{D}_{\rho_1-\frac{\delta_1}{2},\,\rho_2-\frac{4}{3}\delta_2}(G)\,.$$
In addition, the conjugate symplectic map $T_{\gamma}'$ can be
generated by a function $\widetilde{F}$ which is analytic in a
domain containing $\mathcal {D}_{\rho-3\delta}(G)\,.$\\[-5pt]

Now, we state and prove an iterative lemma, which is
then used in the proof of Normal Form Lemma.

\textbf{Lemma 2.8 (Iterative Lemma)} Consider
$F(\hat{x},y)=F_0(\hat{x})+Z(\hat{x},y)+R(\hat{x},y)$ real analytic
in $\mathcal {D}_{\rho}(G)$ and $\mathcal {V}_{\rho_1}(G)$ is
$\alpha,K$-nonresonance modulo $\mathcal {M}$, $Z\in \mathcal
{R}(\mathcal {M},K)$. Let the symplectic map
$T_{\gamma}:(x,y)\rightarrow (\hat{x},\hat{y})$ be given by the
generating function $F(\hat{x},y)$, and let
$\widetilde{\omega}(x)=\partial F_0(x)$ satisfying the condition
(\ref{eq:3}) where $\beta\leq \alpha$. Assume
\begin{equation}\label{eq:7}
||DZ||_{G,\,\rho,\,c}+||DR||_{G,\,\rho,\,c}\leq\frac{\alpha\hat{\delta}_c}{4A}
\quad \mbox{with}\quad 3\delta\leq\rho\,, ~\rho_1\leq
\min\Big(\frac{\delta_2}{6\beta},\,1\Big)
\end{equation}
where $A=1+\frac{\beta c}{\alpha}e^{K\beta\rho_1}$\,. Then, there
exists a real analytic canonical transformation $\Phi
:(a,\varphi)\rightarrow (x,y)$ such that the conjugate symplectic
map $T_{\gamma}'=\Phi^{-1}\circ T_{\gamma}\circ\Phi$ is
generated by
$\widetilde{F}(\hat{a},\varphi)=F_0(\hat{a})+\widetilde{Z}(\hat{a},\varphi)+\widetilde{R}(\hat{a},\varphi)$
with $\widetilde{Z}\in \mathcal {R}(\mathcal {M},K)$ and
$\widetilde{F}$ is analytic in a domain containing $\mathcal
{D}_{\rho-3\delta}(G)$. In addition, ones have
\begin{description}
  \item[(a)] $||D\widetilde{Z}||_{G,\,\rho-3\delta,\,c}\leq
||DZ||_{G,\,\rho,\,c}+||DR||_{G,\,\rho,\,c}\,.$
  \item[(b)] $||D\widetilde{R}||_{G,\,\rho-3\delta,\,c}\leq
\Big\{\frac{2}{\delta_1}\big(1+\frac{A}{\alpha}\big)\big[\triangle_1(||DZ||_{G,\,\rho,\,c}+||DR||_{G,\,\rho,\,c})+e^{-K\delta_2}\big]\cdot
||DR||_{G,\,\rho,\,c}+\triangle_2(||DZ||_{G,\,\rho,\,c}+||DR||_{G,\,\rho,\,c})+e^{-K\delta_2}+\frac{2\beta
AC}{\alpha}\Big\}\cdot||DR||_{G,\,\rho,\,c}$\,, where
\begin{align*}
\triangle_1&=\frac{2A}{\alpha\delta_1}+\frac{2}{\delta_1}+\frac{A}{e\alpha\delta_2}\Big(\frac{\beta
A}{\alpha}+\beta+\frac{2}{c}\Big)\,,\\
\triangle_2&=\frac{n+2}{\delta_1}\cdot\frac{A}{\alpha}+\frac{n+1}{\delta_1}+\frac{nA}{e\alpha\delta_2}\Big(\frac{\beta
A}{\alpha}+\beta+\frac{1}{c}\Big)\,.
\end{align*}
  \item[(c)] $|P_{x}\Phi-id|_{G,\,\rho-3\delta}\leq
\frac{A}{\alpha}||DR||_{G,\,\rho,\,c}\,.$
\end{description}

\textbf{Proof.} We define the transformation
$\Phi:(a,\varphi)\rightarrow (x,y)$ implicitly with the help of a
undetermined generating function $\chi$ by
$$x=a+\partial_2\chi(a,y)\,,\qquad y=\varphi-\partial_1 \chi(a,y)\,.$$
Because of Lemma 2.7, we have the generating function
$\widetilde{F}(\hat{a},\varphi)$ of the conjugate symplectic map
$T_{\gamma}'$ and $\widetilde{F}$ is analytic in domain containing
$\mathcal {D}_{\rho-3\delta}(G)$.

We choose $\chi$ satisfying the linear functional equation:
$$R^{\leq K}+\chi(a,y+\widetilde{\omega}(a))-\chi(a,y)=g(a,\varphi)\qquad \mbox{for some}\quad g\in \mathcal {R}(\mathcal {M},K)\,,$$
where $R^{\leq K}$ denotes the terms restricted to $|k|_1\leq K$ in
the Fourier expansion of $R$. By means of the Fourier expansions of
$R$, $\chi$ and $g$, ones have the solutions:
\begin{equation*}\chi_k(a)=
\begin{cases}\frac{R_k(a)}{1-e^{ik\cdot\widetilde{\omega}(a)}} & ~k\notin \mathcal
{M},~|k|_1\leq K.\\0 &~\text{else.}
\end{cases}
\end{equation*}
and
\begin{equation*}g_k(a)=
\begin{cases}R_k(a) & ~k\in \mathcal
{M},~|k|_1\leq K.\\0 &~\text{else.}
\end{cases}
\end{equation*}\\[3pt]
where
$$\chi=\sum_{k\in \mathbb{Z}^n}\chi_k(a)e^{ik\cdot y}\,,\quad g=\sum_{k\in \mathbb{Z}^n}\chi_k(a)e^{ik\cdot y}\,.$$
Thus, $ \chi(a,y)=\sum\limits_{\begin{subarray}{l} k\notin \mathcal
{M}\\|k|_1 \leq K
\end{subarray}}\frac{R_k(a)}{1-e^{ik\cdot\widetilde{\omega}(a)}}\cdot
e^{ik\cdot y}\,, $ and
\begin{equation*}
||\partial_2\chi||_{G,\,\rho,\,1}=\sum_{\begin{subarray}{l} k\notin
\mathcal {M}\\|k|_1 \leq K
\end{subarray}}\Big|\frac{kR_k(a)}{1-e^{ik\cdot\widetilde{\omega}(a)}}\Big|_{G,\,\rho_1,\,1}\cdot
e^{|k|_1\cdot\rho_2}\leq \frac{1}{\alpha}||\partial_2
R||_{G,\,\rho,\,1}\,,
\end{equation*}
where in the last inequality we have used the nonresonance
condition. Because
$$\partial \chi_k(a)=\frac{\partial R_k(a)}{1-e^{ik\cdot\widetilde{\omega}(a)}}+
\frac{e^{ik\cdot\widetilde{\omega}(a)}[\partial_2R]_k\partial
\widetilde{\omega}}{(1-e^{ik\cdot\widetilde{\omega}(a)})^2}\,,\qquad
\mbox{for}~k\notin \mathcal {M},~|k|_1\leq K\,,$$ where we have used
that $[\partial_2R]_k=iR_k(a)k$ (differentiating the Fourier
expansion of $R$). As before,
$|\texttt{Im}\,\widetilde{\omega}(a)|_{G,\,\rho_1,\,\infty}\leq
\beta\rho_1$. Thus
$$|\partial\chi_k(a)|_{G,\,\rho_1,\,\infty}\leq \frac{1}{\alpha}|\partial R_k(a)|_{G,\,\rho_1,\,\infty}
+\frac{\beta}{\alpha^2}e^{|k|_1\beta\rho_1}\cdot|[\partial_2
R]_k|_{G,\,\rho_1}\,.$$ Moreover,
$$||\partial_1\chi||_{G,\,\rho,\,\infty}\leq \frac{1}{\alpha}||\partial_1 R|_{G,\,\rho,\,\infty}
+\frac{\beta}{\alpha^2}e^{K\beta\rho_1}\cdot||\partial_2
R||_{G,\,\rho,\,1}\,.$$ Therefore, we have
\begin{equation}\label{eq:8}
\begin{split}
||D\chi||_{G,\,\rho,\,c}&\leq
\big(\frac{1}{\alpha}+\frac{\beta c}{\alpha^2}e^{K\beta\rho_1}\big)\cdot||D
R||_{G,\,\rho,\,c} \\
 &\triangleq \frac{A}{\alpha}||DR||_{G,\,\rho,\,c}
 \end{split}
 \end{equation}
with $A=1+\frac{\beta c}{\alpha}e^{K\beta\rho_1}$.\\[-5pt]

Let $B=R_0^{-1}\circ\Phi^{-1}\circ T_{\gamma}\circ\Phi
:(a,\varphi)\rightarrow (\hat{a},\hat{\varphi})$ where $R_0$ is a
integrable rotation on $G\times T^n$ with frequency map
$\widetilde{\omega}$, i.e.
$R_0(a,\varphi)=(a,\varphi+\widetilde{\omega}(a))$, and $B$ can be
expressed implicitly as follows:
\begin{equation}\label{eq:31}
\left\{ \begin{aligned}
         \hat{a} &=a+\partial_2\chi(a,y)-\partial_2\chi(\hat{a},\hat{y})-\partial_2(Z+R)(\hat{x},y)  \\
                 \hat{\varphi}&=\varphi-\partial_1\chi(a,y)+\partial_1\chi(\hat{a},\hat{y})+\partial_1(Z+R)(\hat{x},y)-\widetilde{\omega}(\hat{a})+\widetilde{\omega}(\hat{x}).
                          \end{aligned} \right.
                          \end{equation}

On the other hand, let
$\widetilde{Z}(\hat{a},\varphi)=Z(\hat{a},\varphi)+\mathcal
{P}_{\mathcal {M}}\mathcal {T}_{K}R(\hat{a},\varphi)$, where
$\mathcal {P}_{\mathcal {M}}\mathcal {T}_{K}R$ denotes the terms
restricted to $k\in \mathcal {M}$ and $|k|_1\leq K$ in the Fourier
expansion of $R$, and assume that the map $\Phi^{-1}\circ
T_{\gamma}\circ\Phi$ has the form:
\begin{equation*}
\left\{ \begin{aligned}
         \hat{a} &= a-\partial_2 \widetilde{Z}(\hat{a},\varphi)-\partial_2 \widetilde{R}(\hat{a},\varphi) \\
                  \hat{\varphi}&=\varphi+\widetilde{\omega}(\hat{a})+\partial_1 \widetilde{Z}(\hat{a},\varphi)+\partial_1
                  \widetilde{R}(\hat{a},\varphi)\,.
                          \end{aligned} \right.
                          \end{equation*}
Then, $B$ has the form:
\begin{equation}\label{eq:32}
\left\{ \begin{aligned}
         \hat{a} &= a-\partial_2 \widetilde{Z}(\hat{a},\varphi)-\partial_2 \widetilde{R}(\hat{a},\varphi) \\
                  \hat{\varphi}&=\varphi+\partial_1 \widetilde{Z}(\hat{a},\varphi)+\partial_1
                  \widetilde{R}(\hat{a},\varphi)\,.
                          \end{aligned} \right.
                          \end{equation}

Combining (\ref{eq:31}) and (\ref{eq:32}), we have
\begin{align*}
||\partial_2
\widetilde{R}(\hat{a},\varphi)||_{G,\,\rho-3\delta,\,1}&\leq
||\hat{a}-a+\partial_2
\widetilde{Z}(\hat{a},\varphi)||_{G,\,\rho-\delta,\,1} \\
&=||\partial_2\chi(a,y)-\partial_2\chi(\hat{a},\hat{y})-\partial_2(Z+R)(\hat{x},y)+\partial_2
\widetilde{Z}(\hat{a},\varphi)||_{G,\,\rho-\delta,\,1}\\
&\leq I_1+I_2+I_3+I_4+I_5+I_6\,,\\[4pt]
||\partial_1
\widetilde{R}(\hat{a},\varphi)||_{G,\,\rho-3\delta,\,\infty}&\leq
||\hat{\varphi}-\varphi-\partial_1
\widetilde{Z}(\hat{a},\varphi)||_{G,\,\rho-\delta,\,\infty} \\
&=||\partial_1\chi(\hat{a},\hat{y})-\partial_1\chi(a,y)+\partial_1(Z+R)(\hat{x},y)-\widetilde{\omega}(\hat{a})+\widetilde{\omega}(\hat{x})\\
&\quad-\partial_1
\widetilde{Z}(\hat{a},\varphi)||_{G,\,\rho-\delta,\,\infty}\\[3pt]
&\leq J_1+J_2+J_3+J_4+J_5+J_6+J_7\,.
\end{align*}
where
\begin{align*}
I_1&=||\partial_2\chi(\hat{a},\hat{y})-\partial_2\chi(a,\hat{y})||_{G,\,\rho-\delta,\,1}\,,\\
I_2&=||\partial_2\chi(a,\hat{y})-\partial_2\chi(a,y+\widetilde{\omega}(a))||_{G,\,\rho-\delta,\,1}\,,\\
I_3&=||\partial_2\chi(a,y+\widetilde{\omega}(a))-\partial_2\chi(a,y)+\partial_2(R(a,y)-\mathcal
{P}_{\mathcal {M}}\mathcal
{T}_{K}R(a,y))||_{G,\,\rho-\delta,\,1}\,,\\
I_4&=||\partial_2R(\hat{x},y)-\partial_2R(a,y)||_{G,\,\rho-\delta,\,1}\,,\\
I_5&=||\partial_2Z(\hat{x},y)-\partial_2Z(\hat{a},y)||_{G,\,\rho-\delta,\,1}+||\partial_2\mathcal
{P}_{\mathcal {M}}\mathcal {T}_{K}R(a,y)-\partial_2\mathcal
{P}_{\mathcal {M}}\mathcal
{T}_{K}R(\hat{a},y)||_{G,\,\rho-\delta,\,1}\,,\\
I_6&=||\partial_2Z(\hat{a},y)-\partial_2Z(\hat{a},\varphi)||_{G,\,\rho-\delta,\,1}+||\partial_2\mathcal
{P}_{\mathcal {M}}\mathcal {T}_{K}R(\hat{a},y)-\partial_2\mathcal
{P}_{\mathcal {M}}\mathcal
{T}_{K}R(\hat{a},\varphi)||_{G,\,\rho-\delta,\,1}\,,\\
J_1&=||\partial_1\chi(\hat{a},\hat{y})-\partial_1\chi(a,\hat{y})||_{G,\,\rho-\delta,\,\infty}\,,\\
J_2&=||\partial_1\chi(a,\hat{y})-\partial_1\chi(a,y+\widetilde{\omega}(a))||_{G,\,\rho-\delta,\,\infty}\,,\\
J_3&=||\partial_1\chi(a,y+\widetilde{\omega}(a))-\partial_1\chi(a,y)+\partial
\widetilde{\omega}(a)\partial_2\chi(a,y+\widetilde{\omega}(a))\\
&\quad+\partial_1(R(a,y)-\mathcal {P}_{\mathcal {M}}\mathcal
{T}_{K}R(a,y))||_{G,\,\rho-\delta,\,\infty}\,,\\
J_4&=||\partial_1R(\hat{x},y)-\partial_1R(a,y)||_{G,\,\rho-\delta,\,\infty}\,,\\
J_5&=||\partial_1Z(\hat{x},y)-\partial_1Z(\hat{a},y)||_{G,\,\rho-\delta,\,\infty}+||\partial_1\mathcal
{P}_{\mathcal {M}}\mathcal {T}_{K}(R(a,y)-R(\hat{a},y))||_{G,\,\rho-\delta,\,\infty}\,,\\
J_6&=||\partial_1Z(\hat{a},y)-\partial_1Z(\hat{a},\varphi)||_{G,\,\rho-\delta,\,\infty}+||\partial_1\mathcal
{P}_{\mathcal {M}}\mathcal {T}_{K}(R(\hat{a},y)-R(\hat{a},\varphi))||_{G,\,\rho-\delta,\,\infty}\,,\\
J_7&=||\widetilde{\omega}(\hat{x})-\widetilde{\omega}(\hat{a})||_{G,\,\rho-\delta,\,\infty}+||\partial
\widetilde{\omega}(a)\cdot\partial_2\chi(a,y+\widetilde{\omega}(a))||_{G,\,\rho-\delta,\,\infty}\,.
\end{align*}
Note that all the concerned variables are in $\mathcal
{D}_{\rho-\delta}(G)$ and we are able to get the following estimates
by middle value theorem and Cauchy estimates.
\begin{align*}
I_1&\leq
\frac{1}{\delta_1}||\partial_2\chi||_{G,\,\rho,\,1}\cdot||\hat{a}-a||_{G,\,\rho-\delta,\,1}\\
&\leq\frac{1}{\delta_1}||D\chi||_{G,\,\rho,\,c}\cdot(||\hat{a}-a+\partial_2\widetilde{Z}(\hat{a},\varphi)||_{G,\,\rho-\delta,\,1}+
||\partial_2\widetilde{Z}(\hat{a},\varphi)||_{G,\,\rho-\delta,\,1})\\
&\leq\frac{1}{\delta_1}\cdot\frac{A}{\alpha}||DR||_{G,\,\rho,\,c}\cdot(||\hat{a}-a+\partial_2\widetilde{Z}(\hat{a},\varphi)||_{G,\,\rho-\delta,\,1}+
||\partial_2\widetilde{Z}(\hat{a},\varphi)||_{G,\,\rho-\delta,\,1})\,,\\
\intertext{where in the last inequality, we have used (\ref{eq:8}).}
I_2&\leq
\frac{1}{e\delta_2}||\partial_2\chi||_{G,\,\rho,\,1}\cdot||\hat{y}-y-\widetilde{\omega}(a)||_{G,\,\rho-\delta,\,\infty}\\
&\leq \frac{1}{e\delta_2}||D\chi||_{G,\,\rho,\,c}\cdot||\hat{y}-y-\widetilde{\omega}(a)||_{G,\,\rho-\delta,\,\infty}\,.
\end{align*}
Note that
\begin{align*}
||\hat{y}-y-\widetilde{\omega}(a)||_{G,\,\rho-\delta,\,\infty}&=||\widetilde{\omega}(\hat{x})-\widetilde{\omega}(a)+
\partial_1(Z+R)(\hat{x},y)||_{G,\,\rho-\delta,\,\infty}\\
&\leq\beta\cdot||\hat{x}-a||_{G,\,\rho-\delta,\,1}+||\partial_1(Z+R)(\hat{x},y)||_{G,\,\rho-\delta,\,\infty}\,,
\end{align*}
and $\hat{x}=a+\partial_2\chi(a,y)-\partial_2(Z+R)(\hat{x},y)$, so
$$||\hat{x}-a||_{G,\,\rho-\delta,\,1}\leq ||\partial_2\chi||_{G,\,\rho-\delta,\,1}+||\partial_2(Z+R)||_{G,\,\rho-\delta,\,1}.$$
Thus,
$$||\hat{y}-y-\widetilde{\omega}(a)||_{G,\,\rho-\delta,\,\infty}\leq \beta(||D\chi||_{G,\,\rho-\delta,\,c}+||D(Z+R)||_{G,\,\rho-\delta,\,c})
+\frac{1}{c}||D(Z+R)||_{G,\,\rho-\delta,\,c}\,.$$
Therefore,
\begin{align*}
I_2&\leq
\frac{A}{e\alpha\delta_2}||DR||_{G,\,\rho,\,c}\cdot\Big[\frac{\beta A}{\alpha}||DR||_{G,\,\rho,\,c}+(\beta+\frac{1}{c})||D(Z+R)||_{G,\,\rho,\,c}\Big]\,.\\
\intertext{Similarly, we have}
I_3&\leq ||\partial_2 R^{>K}||_{G,\,\rho-\delta,\,1}\leq
||DR^{>K}||_{G,\,\rho-\delta,\,c}\leq
e^{-K\delta_2}||DR||_{G,\,\rho,\,c}\,,\\
I_4&\leq
\frac{1}{\delta_1}||\partial_2R||_{G,\,\rho,\,1}\cdot||\hat{x}-a||_{G,\,\rho-\delta,\,1}\leq \frac{1}{\delta_1}||DR||_{G,\,\rho,\,c}
\cdot\Big(\frac{A}{\alpha}||DR||_{G,\,\rho,\,c}+||D(Z+R)||_{G,\,\rho,\,c}\Big),\\
I_5&\leq
\frac{1}{\delta_1}||\partial_2Z||_{G,\,\rho,\,1}\cdot||\hat{x}-\hat{a}||_{G,\,\rho-\delta,\,1}
+\frac{1}{\delta_1}||\partial_2R||_{G,\,\rho,\,1}\cdot||\hat{a}-a||_{G,\,\rho-\delta,\,1}\\
&\leq \frac{1}{\delta_1}||DZ||_{G,\,\rho,\,c}
\cdot\frac{A}{\alpha}||DR||_{G,\,\rho,\,c}+\frac{1}{\delta_1}||DR||_{G,\,\rho,\,c}
\cdot\big(||\hat{a}-a+\partial_2\widetilde{Z}||_{G,\,\rho-\delta,\,1}+||\partial_2\widetilde{Z}||_{G,\,\rho,\,1}\big)\,,\\
I_6&\leq\frac{1}{e\delta_2}\cdot\frac{A}{c\alpha}||DR||_{G,\,\rho,\,c}\cdot\big(||DZ||_{G,\,\rho,\,c}+||DR||_{G,\,\rho,\,c}\big)\,.
\end{align*}
By means of the condition (\ref{eq:7}) and above estimates, we get
\begin{align*}
\frac{1}{2}||\hat{a}-a+\partial_2\widetilde{Z}(\hat{a},\varphi)||_{G,\,\rho-\delta,\,1}\leq& \bigg[\Big(\frac{2A}{\alpha\delta_1}+\frac{2}{\delta_1}+\frac{A}{e\alpha\delta_2}\big(\beta\frac{A}{\alpha}+\beta+
\frac{2}{c}\big)\Big)\cdot\Big(||DZ||_{G,\,\rho,\,c}\\
&+||DR||_{G,\,\rho,\,c}\Big)+e^{-K\delta_2}\bigg]\cdot||DR||_{G,\,\rho,\,c}\,.
\end{align*}
Thus,
$$||\partial_2\widetilde{R}(\hat{a},\varphi)||_{G,\,\rho-3\delta,\,1}\leq 2\big[\triangle_1\cdot(||DR||_{G,\,\rho,\,c}+||DZ||_{G,\,\rho,\,c})+e^{-K\delta_2}\big]\cdot||DR||_{G,\,\rho,\,c}\,,$$
where $\triangle_1=\frac{2A}{\alpha\delta_1}+\frac{2}{\delta_1}+\frac{A}{e\alpha\delta_2}(\frac{\beta
A}{\alpha}+\beta+\frac{2}{c})\,.$

By estimating $J_k, ~k=1,\ldots,7$ in a similar way to the above and making use of the previous estimates, we obtain
\begin{align*}
J_1&\leq
\frac{1}{\delta_1}||\partial_1\chi||_{G,\,\rho,\,\infty}\cdot||\hat{a}-a||_{G,\,\rho-\delta,\,1}\\
&\leq\frac{1}{c\delta_1}\cdot\frac{A}{\alpha}||DR||_{G,\,\rho,\,c}\cdot\big(||\hat{a}-a+\partial_2\widetilde{Z}||_{G,\,\rho-\delta,\,1}+
||\partial_2\widetilde{Z}||_{G,\,\rho-\delta,\,1}\big)\,,\\
J_2&\leq
\frac{1}{e\delta_2}||\partial_1\chi||_{G,\,\rho}\cdot||\hat{y}-y-\widetilde{\omega}(a)||_{G,\,\rho-\delta,\,\infty}\\
&\leq
\frac{n}{e\delta_2}\cdot\frac{A}{c\alpha}||DR||_{G,\,\rho,\,c}\cdot\Big(\frac{\beta A}{\alpha}||DR||_{G,\,\rho,\,c}+\big(\beta+\frac{1}{c}\big)||D(Z+R)||_{G,\,\rho,\,c}\Big)\,,\\
J_3&\leq ||\partial_1 R^{>K}||_{G,\,\rho-\delta,\,\infty}\leq\frac{1}{c}
||DR^{>K}||_{G,\,\rho-\delta,\,c}\leq\frac{1}{c}
e^{-K\delta_2}||DR||_{G,\,\rho,\,c}\,,\\
J_4&\leq
\frac{n}{c\delta_1}||DR||_{G,\,\rho,\,c}
\cdot\Big(\frac{A}{\alpha}||DR||_{G,\,\rho,\,c}+||D(Z+R)||_{G,\,\rho,\,c}\Big),\\
J_5
&\leq \frac{1}{c\delta_1}||DZ||_{G,\,\rho,\,c}
\cdot\frac{A}{\alpha}||DR||_{G,\,\rho,\,c}+\frac{1}{c\delta_1}||DR||_{G,\,\rho,\,c}
\cdot\big(||\hat{a}-a+\partial_2\widetilde{Z}||_{G,\,\rho-\delta,\,1}+||\partial_2\widetilde{Z}||_{G,\,\rho,\,1}\big)\,,\\
J_6&\leq\frac{1}{\delta_1}\cdot\frac{A}{c\alpha}||DR||_{G,\,\rho,\,c}\cdot\big(||DZ||_{G,\,\rho,\,c}+||DR||_{G,\,\rho,\,c}\big)\,,\\
J_7&\leq 2\beta||\partial_2\chi||_{G,\,\rho-\delta,\,1}\leq 2\beta\cdot\frac{A}{\alpha}||DR||_{G,\,\rho,\,c}\,.
\end{align*}
Combining the above estimates, we get
\begin{align*}
c||\partial_1
\widetilde{R}(\hat{a},\varphi)||_{G,\,\rho-3\delta,\,\infty}\leq & \Big[\frac{1}{\delta_1}\big(1+\frac{A}{\alpha}\big)\cdot||\hat{a}-a+\partial_2\widetilde{Z}||_{G,\,\rho-\delta,\,1}+
\triangle_2\cdot\big(||DZ||_{G,\,\rho,\,c}+||DR||_{G,\,\rho,\,c}\big)\\
&+e^{-K\delta_2}+\frac{2\beta Ac}{\alpha}\Big]\cdot||DR||_{G,\,\rho,\,c}\,,
\end{align*}
where $$\triangle_2=\frac{n+2}{\delta_1}\cdot\frac{A}{\alpha}+\frac{n+1}{\delta_1}+\frac{nA}{e\alpha\delta_2}\Big(\frac{\beta
A}{\alpha}+\beta+\frac{1}{c}\Big)\,.$$
Because $\delta_1\leq 1$, we have
\begin{align*}
||DR||_{G,\,\rho-3\delta,\,c}&\leq \max\big(||\partial_2
\widetilde{R}||_{G,\,\rho-3\delta,\,1},~c||\partial_1
\widetilde{R}||_{G,\,\rho-3\delta,\,\infty}\big)\\
&\leq\Big\{\frac{2}{\delta_1}\big(1+\frac{A}{\alpha}\big)\big[\triangle_1(||DZ||_{G,\,\rho,\,c}+||DR||_{G,\,\rho,\,c})+e^{-K\delta_2}\big]\cdot
||DR||_{G,\,\rho,\,c}\\&\quad+\triangle_2(||DZ||_{G,\,\rho,\,c}+||DR||_{G,\,\rho,\,c})+e^{-K\delta_2}+\frac{2\beta
Ac}{\alpha}\Big\}\cdot||DR||_{G,\,\rho,\,c}\,.
\end{align*}

Finally,
$$|P_{x}\Phi-id|_{G,\,\rho-3\delta}\leq||\partial_2\chi||_{G,\,\rho-\delta,\,1}\leq ||D\chi||_{G,\,\rho,\,c}\leq
\frac{A}{\alpha}||DR||_{G,\,\rho,\,c}\,.$$
\qed

\subsection{Proof of Normal Form Lemma}

As the Hamiltonian case, we shall construct a series of symplectic transformations $\Phi^{(i)}$, each of which reduces the norm of remainder by factor $\frac{1}{e}$. After applying Iterative Lemma $N$ times, we can get an exponentially small remainder by choosing $N=N(K)$ adequately.

Let $N\geq 1$ be an integer to be chosen below. Denoting $\rho^{(i)}=\rho-3i\delta$, with $\delta=\frac{\rho}{6N}$. Obviously, $\rho^{(i)}=\rho^{(i-1)}-3\delta$.

Next we apply Iterative Lemma $N$ times, and obtain a series of symplectic transformations $\Phi^{(i)}:\mathcal {D}_{\rho^{(i)}}(G)\rightarrow
\mathcal {D}_{\rho^{(i-1)}}(G)$ for $1\leq i\leq N$. Let $\Psi^{(i)}=\Phi^{(i)}\circ\Phi^{(i-1)}\circ\cdots\circ\Phi^{(1)}$, $T_{\gamma}^{(i)}={\Psi^{(i)}}^{-1}\circ T_{\gamma}\circ\Psi^{(i)}$ and the generating functions $F^{(i)}=F_0+Z^{(i)}+R^{(i)}$ of the symplectic map $T_{\gamma}^{(i)}$ with $Z^{(i)}\in \mathcal {R}(\mathcal {M},K)$.

Now, we are going to show that if $\frac{K\rho_2}{6N}\geq 2$ then the claims below are true for $1\leq i\leq N$.
\begin{description}
  \item[(a)] $||DZ^{(i)}||_{G,\,\rho^{(i)},\,c}\leq \sum\limits_{j=0}^{i-1}||DR^{(j)}||_{G,\,\rho^{(j)},\,c}\,.$
  \item[(b)] $||DR^{(i)}||_{G,\,\rho^{(i)},\,c}\leq {\Large\frac{1}{e}}||DR^{(i-1)}||_{G,\,\rho^{(i-1)},\,c}\,.$
\end{description}

The proof is done by induction.

Setting $R^{(0)}=f$ and $Z^{(0)}=0$, and we choose the parameter $c=\frac{\rho_1}{\rho_2}=\frac{\delta_1}{\delta_2}$ so that $\hat{\delta}_c=\delta_1=c\delta_2$.
Due to $\frac{1}{K}\leq \frac{\rho_2}{12N}$ and $\widetilde{C}\geq 30$, we have
$$||Df||_{G,\,\rho,\,c}\leq\frac{\alpha\rho_1}{\widetilde{C}AK\rho_2}\leq\frac{\alpha\rho_1}{12\widetilde{C}AN}=
\frac{\alpha\delta_1}{2\widetilde{C}A}\leq\frac{\alpha\hat{\delta}_c}{60A}\,.$$
Note that $$\rho_1\leq \frac{\alpha}{4K\beta}\leq \frac{\alpha\rho_2}{48N\beta}=\frac{\alpha\delta_2}{8\beta}\leq \frac{\delta_2}{6\beta}\,,$$
Thus, Iterative Lemma can be applied with $\frac{\rho}{6N}$ instead of $\delta$. Due to $\rho_1\leq
\frac{\alpha}{4K\beta}$, we have
\begin{equation}\label{eq:9}
A=1+\frac{\beta c}{\alpha}e^{K\beta\rho_1}\leq 1+\frac{2\beta c}{\alpha}\leq 1+\frac{1}{2K\rho_2}\leq 2 \,,
\end{equation}
and
\begin{equation}\label{eq:10}
\frac{2\beta Ac}{\alpha}\leq \frac{1}{K\rho_2}\leq \frac{1}{12}\,.
\end{equation}
By means of (\ref{eq:5}) and (\ref{eq:9}), we get
\begin{align*}
||DR^{(1)}||_{G,\,\rho^{(1)},\,c}&\leq \Big(\frac{1}{e\widetilde{C}}+\frac{n+2}{2\widetilde{C}}+\frac{n+1}{2\widetilde{C}}+\frac{n}{24e\widetilde{C}}
+\frac{1}{e^2}+\frac{1}{12}\Big)\cdot||Df||_{G,\,\rho,\,c}\\
&\leq \frac{1}{e}||Df||_{G,\,\rho,\,c}\,.
\end{align*}
The claim (a) is obviously true for $i=1$.

For $1<i\leq N$, note that
\begin{align*}
||DR^{(i-1)}||_{G,\,\rho^{(i-1)},\,c}&\leq \frac{||DR^{(0)}||_{G,\,\rho,\,c}}{e^{i-1}}\leq ||Df||_{G,\,\rho,\,c}\leq \frac{\alpha\hat{\delta}_c}{60A}\,,\\
\intertext{and}
||DZ^{(i-1)}||_{G,\,\rho^{(i-1)},\,c}&\leq \sum\limits_{j=0}^{i-2}||DR^{(j)}||_{G,\,\rho^{(j)},\,c}\leq 2||DR^{(0)}||_{G,\,\rho,\,c}
\leq \frac{\alpha\hat{\delta}_c}{30A}\,.
\end{align*}

Thus, Iterative Lemma can be applied with $\frac{\rho}{6N}$ instead of $\delta$. The claim (a) is easy to prove, and the claim (b) can be proved by the following estimates:
\begin{align*}
||DR^{(i)}||_{G,\,\rho^{(i)},\,c}&\leq
\bigg\{\frac{2}{\delta_1}\big(1+\frac{A}{\alpha}\big)\Big[\triangle_1\big(||DZ^{(i-1)}||_{G,\,\rho^{(i-1)},\,c}+||DR^{(i-1)}||_{G,\,\rho^{(i-1)},\,c}\big)
+e^{-K\delta_2}\Big]\cdot\\
&\qquad ||DR^{(i-1)}||_{G,\,\rho^{(i-1)},\,c}
+\triangle_2\big(||DZ^{(i-1)}||_{G,\,\rho^{(i-1)},\,c}+||DR^{(i-1)}||_{G,\,\rho^{(i-1)},\,c}\big)\\
&\qquad +e^{-K\delta_2}+\frac{2\beta
AC}{\alpha}\bigg\}\cdot||DR^{(i-1)}||_{G,\,\rho^{(i-1)},\,c}\\
&\leq
\bigg\{\frac{2}{\delta_1}\big(1+\frac{A}{\alpha}\big)\Big(\triangle_1\cdot 2||Df||_{G,\,\rho,\,c}+e^{-\frac{K\rho_2}{6N}}\Big)\frac{||Df||_{G,\,\rho,\,c}}{e}+2\triangle_2||Df||_{G,\,\rho,\,c}\\
&\qquad +e^{-\frac{K\rho_2}{6N}}+\frac{1}{12}\bigg\}\cdot||DR^{(i-1)}||_{G,\,\rho^{(i-1)},\,c}\\
&\leq \frac{1}{e}||DR^{(i-1)}||_{G,\,\rho^{(i-1)},\,c}\,.
\end{align*}

Now we may choose $N=\big[\frac{K\rho_2}{12}\big]$, the integer part of $\frac{K\rho_2}{12}$. After iterating $N$ times, the exponential small remainder is given by
$$||DR||_{G,\,\frac{\rho}{2},\,c}=||DR^{(N)}||_{G,\,\frac{\rho}{2},\,c}\leq \frac{||Df||_{G,\,\rho,\,c}}{e^N}\leq 3e^{-\frac{K\rho_2}{12}}||Df||_{G,\,\rho,\,c}\,.$$

The conclusion 3) is obtained from the fact that
\begin{align*}
|P_{x}\Phi-id|_{G,\,\frac{\rho}{2}}&\leq \sum_{i=1}^{N}|P_{x}\Phi^{(i)}-id|_{G,\,\rho^{i}}\leq \sum_{i=1}^{N}\frac{A}{\alpha}||DR^{(i)}||_{G,\,\rho^{(i)},\,c}\leq \frac{2A}{\alpha}||Df||_{G,\,\rho,\,c}\\
&\leq\frac{\rho_1}{2^8}\,,
\end{align*}
where the last inequality is a consequence of (\ref{eq:5}) and $\widetilde{C}\geq 48$. Here we remark that if $K\rho_2\leq 12$, all results are obvious if we take $\Phi$ as the identity map.

\section{The Geometry of Resonances}

In this section, we concern the covering of the whole action space $\mathcal {G}_{\gamma}$ by a family of resonant blocks associated to different lattice $\mathcal {M}$. For the symplectic map, the original geometric construction in \cite{Benettin} requires some modifications (see also \cite{Guzzo}). Here in addition to $K$, our geometric construction will be characterized  by $2n$ positive parameters $0<\alpha_0=\alpha_1\leq\alpha_2\leq\cdots\leq\alpha_n<1$ and $\delta_1, \delta_2,\cdots,\delta_n$. More precisely, for any choice of these parameters, for each $K$-lattice $\mathcal {M}\subseteq\mathbb{Z}^n$ with dim$\mathcal {M}=r$, and some $l=(l_1,\ldots,l_r)\in\mathbb{Z}^r$, we define:

\textbf{i) Resonant manifold}
$$\mathcal {R}_{\mathcal {M}}^l=\big\{x\in\mathcal {G}_{\gamma};~k^{(j)}\cdot\widetilde{\omega}(x)+2\pi l_j=0,~j=1,\ldots,r \big\}.$$
where $k^{(1)},\ldots,k^{(r)}$ is a $K$-basis.

\textbf{ii) Resonant zone}
$$\mathcal {Z}_{\mathcal {M}}^l=\big\{x\in\mathcal {G}_{\gamma};~\big|k^{(j)}\cdot\widetilde{\omega}(x)+2\pi l_j\big|\leq \alpha_r,~j=1,\ldots,r \big\}.$$
where $k^{(1)},\ldots,k^{(r)}$ is a $K$-basis.

Note that, for $r=0$, $K$-lattice $\mathcal {M}$ is trivial, and $\mathcal {Z}_0$ is defined as coinciding with the whole action space $\mathcal {G}_{\gamma}$. Denoting by $\mathcal {Z}_r^*$, $1\leq r\leq n$, the union of all resonant zones with the same dimension $r$, i.e.
$$\mathcal {Z}_r^*=\bigcup_{l}\bigcup_{\text{dim}(\mathcal {M})=r}\mathcal {Z}_{\mathcal {M}}^l\,,$$
and set $\mathcal {Z}_{n+1}^*=\emptyset$.

\textbf{iii) Resonant block}
$$\mathcal {B}_{\mathcal {M}}^l=\mathcal {Z}_{\mathcal {M}}^l\backslash\mathcal {Z}_{r+1}^*\,.$$
Especially, $\mathcal {B}_0=\mathcal {Z}_{0}\backslash\mathcal {Z}_{1}^*$. The dimension $r$ of $\mathcal {M}$ will also be called the multiplicity of the corresponding resonant manifold, zone or block.

\textbf{iv) Cylinder}

First, let $\Pi_{\mathcal {M}}(x)$ be the hyperplane through $x$ parallel to $\mathcal {M}$ with the same dimensionality, and denote its $\delta_r$ neighborhood by $\Pi_{\mathcal {M},\,\delta_r}(x)$, i.e.
$$\Pi_{\mathcal {M},\,\delta_r}(x)=\big\{\tilde{x}\in \mathbb{R}^n;~\text{dist}(\tilde{x},\Pi_{\mathcal {M}}(x))\leq\delta_r\big\}.$$
Then, for $x\in\mathcal {B}_{\mathcal {M}}^l$, the cylinder is defined by
$$\mathcal {C}_{\mathcal {M},\,\delta_r}^l(x)=\Pi_{\mathcal {M},\,\delta_r}(x)\bigcap\mathcal {Z}_{\mathcal {M}}^l\,,$$
ËüµÄ»ùΪ~$\Pi_{\mathcal {M},\,\delta_r}(x)\bigcap\partial\mathcal {Z}_{\mathcal {M}}$\,.

\textbf{v) Extended resonant block}
$$\mathcal {B}_{\mathcal {M},\,\delta_r}^l=\bigcup_{x\in\mathcal {B}_{\mathcal {M}}^l}\mathcal {C}_{\mathcal {M},\,\delta_r}^l(x)\,.$$

\textbf{Remarks:} 1. In the definition of the resonant zone, we make the point that $l$ is bounded because of $|k|_1\leq K$ and the boundedness of $\widetilde{\omega}(x)$ and $\alpha$. In addition, the resonant zones with the same lattice $\mathcal {M}$ don't intersect for different $l$.

2. The resonant blocks $\mathcal {B}_{\mathcal {M}}^l$ constitute a covering of the action space $\mathcal {G}_{\gamma}$\,, that is
$\mathcal {G}_{\gamma}=\bigcup_{l}\bigcup_{\mathcal {M}}\mathcal {B}_{\mathcal {M}}^l\,.$\\[-5pt]

Now, we shall prove some properties of the geometry construction.

\textbf{Proposition 3.1} (i) For any $x\in\mathcal {B}_{\mathcal {M}}^l$ with dim$(\mathcal {M})=r\in\{0,1,\ldots,n-1\}$, if $k\in\mathbb{Z}_{K}^n\backslash\mathcal {M}$, then $|k\cdot\widetilde{\omega}(x)+2\pi l_0|> \alpha_{r+1}$ for any $l_0\in \mathbb{Z}$. In particular, for any $x\in\mathcal {B}_0$ it is $|k\cdot\widetilde{\omega}(x)+2\pi l_0|> \alpha_{1}$ for any $k\in\mathbb{Z}_{K}^n\backslash\{0\},\,l_0\in \mathbb{Z}$.\\
(ii)
$$\bigcup_{l}\bigcup_{\text{dim}(\mathcal {M})=r}\mathcal {B}_{\mathcal {M}}^l=\mathcal {Z}_{r}^*\backslash\mathcal {Z}_{r+1}^*\,.$$
(iii)
$$\mathcal {G}_{\gamma}\backslash\mathcal {Z}_{r+1}^*=\bigcup_{l}\bigcup_{\text{dim}(\mathcal {M})\leq r}\mathcal {B}_{\mathcal {M}}^l\,.$$
(iv) If $\mathcal {M}$ is $r$-dimensional $K$-lattice, $1\leq r\leq n$, then for any $x\in\mathcal {B}_{\mathcal {M}}^l$,
$$\text{diam}(\mathcal {C}_{\mathcal {M},\,\delta_r}^{l}(x))\leq \frac{4}{\mu}\delta_r+\frac{2rK^{r-1}}{\beta\mu}\alpha_r\,.$$

\textbf{Proof.} In order to prove the statement (i) by contradiction, assume that there exist $l_0\in\mathbb{Z}$ and $x\in\mathcal {B}_{\mathcal {M}}^l$, such that $|k\cdot\widetilde{\omega}(x)+2\pi l_0|\leq \alpha_{r+1}$ for any $k\in\mathbb{Z}_{K}^n\backslash\mathcal {M}$\,. Note that $\alpha_{r}\leq\alpha_{r+1}$, then we obtain that there exists a $r+1$-dimensional $K$-lattice $\mathcal {M}'$ and $l'=(l,l_0)$ such that $x\in\mathcal {Z}_{\mathcal {M}'}^{l'}$\,. However, $x\notin \mathcal {Z}_{r+1}^*$ in view of the definition of $\mathcal {B}_{\mathcal {M}}^l$\,. We get the contradiction.

For the properties (ii) and (iii), it is readily to obtain from the definitions. In fact,
\begin{align*}
\bigcup_{l}\bigcup_{\text{dim}(\mathcal {M})=r}\mathcal {B}_{\mathcal {M}}^l&=\bigcup_{l}\bigcup_{\text{dim}(\mathcal {M})=r}\big(\mathcal {Z}_{\mathcal {M}}^l\backslash\mathcal {Z}_{r+1}^*\big)\\
&=\Big(\bigcup_{l}\bigcup_{\text{dim}(\mathcal {M})=r}\mathcal {Z}_{\mathcal {M}}^l\Big)\backslash\mathcal {Z}_{r+1}^*\\
&=\mathcal {Z}_{r}^*\backslash\mathcal {Z}_{r+1}^*\,,
\end{align*}
and
\begin{align*}
\bigcup_{l}\bigcup_{\text{dim}(\mathcal {M})\leq r}\mathcal {B}_{\mathcal {M}}^l&=\bigcup_{s=0}^{r}\Big(\bigcup_{l}\bigcup_{\text{dim}(\mathcal {M})=s}\mathcal {B}_{\mathcal {M}}^l\Big)\\
&=\bigcup_{s=0}^{r}\Big(\mathcal {Z}_{s}^*\backslash\mathcal {Z}_{s+1}^*\Big)\\
&=\mathcal {G}_{\gamma}\backslash\mathcal {Z}_{r+1}^*\,.
\end{align*}
Finally, before proving the statement (iv), we need a technique lemma which refers to \cite{Benettin}.

\textbf{Lemma 3.1 \cite{Benettin}} Let $k^{(1)},k^{(2)},\cdots,k^{(r)}$ be linearly independent vectors of $\mathbb{Z}^n$ with $|k^{(i)}|_1\leq K$, and $\omega\in \mathbb{R}^n$ be any linear combination of $k^{(1)},k^{(2)},\cdots,k^{(r)}$ satisfying $|\omega\cdot k^{(i)}|\leq \alpha,~(i=1,\ldots,r)$ then one has
$$|\omega|_2<rK^{r-1}\alpha\,.$$

Let's continue the proof of Proposition 3.1\,(iv). For any $x_1,x_2\in\mathcal {C}_{\mathcal {M},\,\delta_r}^{l}(x)$ with $x\in\mathcal {B}_{\mathcal {M}}^l$, there exist $x_1^*,x_2^*\in\Pi_{\mathcal {M}}(x)$ such that $|x_1^*-x_1|_2\leq\delta_r$ and $|x_2^*-x_2|_2\leq\delta_r$\,. Due to the convexity of the function $F_0$, we conclude that
\begin{align}\label{eq:12}
\beta\mu|x_1-x_2|_2^2&\leq \big|(\widetilde{\omega}(x_1)-\widetilde{\omega}(x_2))\cdot(x_1-x_2)\big| \nonumber\\
&\leq\big|(\widetilde{\omega}(x_1)-\widetilde{\omega}(x_2))\cdot(x_1-x_1^*)\big|+\big|(\widetilde{\omega}(x_1)-\widetilde{\omega}(x_2))
\cdot(x_1^*-x_2^*)\big| \nonumber\\
&\quad+\big|(\widetilde{\omega}(x_1)-\widetilde{\omega}(x_2))\cdot(x_2^*-x_2)\big|\nonumber\\
&\leq 2\beta|x_1-x_2|_2\cdot\delta_r+\big|(\widetilde{\omega}(x_1)-\widetilde{\omega}(x_2))
\cdot(x_1^*-x_2^*)\big|\,.
\end{align}
Note that $x_1^*-x_2^*$ parallel to $\mathcal {M}$, so
\begin{equation}\label{eq:13}
\big|(\widetilde{\omega}(x_1)-\widetilde{\omega}(x_2))
\cdot(x_1^*-x_2^*)\big|=\big|\mathcal {P}_{\mathcal {M}}(\widetilde{\omega}(x_1)-\widetilde{\omega}(x_2))
\cdot(x_1^*-x_2^*)\big|
\end{equation}
where $\mathcal {P}_{\mathcal {M}}$ denote the projection of a vector onto $\mathcal {M}$. Moreover,
\begin{align}\label{eq:14}
&\big|\mathcal {P}_{\mathcal {M}}(\widetilde{\omega}(x_1)-\widetilde{\omega}(x_2))
\cdot(x_1^*-x_2^*)\big| \nonumber\\
\leq&\big|\mathcal {P}_{\mathcal {M}}(\widetilde{\omega}(x_1)-\widetilde{\omega}(x_2))
\cdot(x_1^*-x_1)\big|+\big|\mathcal {P}_{\mathcal {M}}(\widetilde{\omega}(x_1)-\widetilde{\omega}(x_2))
\cdot(x_1-x_2)\big|  \nonumber\\
&+\big|\mathcal {P}_{\mathcal {M}}(\widetilde{\omega}(x_1)-\widetilde{\omega}(x_2))
\cdot(x_2-x_2^*)\big| \nonumber\\
\leq& 2\beta|x_1-x_2|_2\cdot\delta_r+\big|\mathcal {P}_{\mathcal {M}}(\widetilde{\omega}(x_1)-\widetilde{\omega}(x_2))
\cdot(x_1-x_2)\big|\,.
\end{align}

Let $k^{(1)},k^{(2)},\cdots,k^{(r)}\in\mathbb{Z}^n$ be the $K$-basis of $\mathcal {M}$. Since $x_1,x_2\in\mathcal {Z}_{\mathcal {M}}^l$, for each $i\in\{1,2,\ldots,r\}$, there exists $l_i\in\mathbb{Z}$ such that $\big|k^{(i)}\cdot\widetilde{\omega}(x_1)+2\pi l_i\big|\leq \alpha_r$ and $\big|k^{(i)}\cdot\widetilde{\omega}(x_2)+2\pi l_i\big|\leq \alpha_r$\,. Thus, we have
\begin{align*}
&\big|\mathcal {P}_{\mathcal {M}}(\widetilde{\omega}(x_1)-\widetilde{\omega}(x_2))
\cdot k^{(i)}\big|=\big|(\widetilde{\omega}(x_1)-\widetilde{\omega}(x_2))
\cdot k^{(i)}\big|\\
=&\big|k^{(i)}\cdot\widetilde{\omega}(x_1)+2\pi l_i-k^{(i)}\cdot\widetilde{\omega}(x_2)-2\pi l_i
\big|\leq 2\alpha_r\,.
\end{align*}
From Lemma 3.1 it follows:
\begin{equation}\label{eq:15}
\big|\mathcal {P}_{\mathcal {M}}(\widetilde{\omega}(x_1)-\widetilde{\omega}(x_2))\big|_2\leq 2rK^{r-1}\alpha_r\,.
\end{equation}
Combining (\ref{eq:12}) -- (\ref{eq:15}), we get
$$\beta\mu|x_1-x_2|_2^2\leq 4\beta|x_1-x_2|_2\cdot\delta_r+2rK^{r-1}\alpha_r\cdot|x_1-x_2|_2\,.$$
That is
$$|x_1-x_2|_2\leq \frac{4}{\mu}\delta_r+\frac{2rK^{r-1}}{\beta\mu}\alpha_r\,.$$
\qed

\section{The Proof of Theorem 1}

In this section, we will complete the proof of Theorem 1. Now, we
have to make a choice for free parameters
$K,\alpha_1,\cdots,\alpha_n$ and $\delta_1,\cdots,\delta_n$ in order
to satisfy two important properties which are crucial points in the
following proof. First, there is no intersection among the extended
blocks with the same dimensional lattices. More precisely, $\mathcal
{V}_{\delta_r}(\mathcal {B}_{\mathcal
{M},\,\delta_r}^l)\bigcap\mathcal {Z}_{\mathcal
{M}'}^{l'}=\emptyset$ for any $\mathcal {M}'\neq\mathcal {M}$ with
$\text{dim}(\mathcal {M})=\text{dim}(\mathcal {M}')=r$, and
$l'\in\mathbb{Z}^r$, that is called the condition of nonoverlapping
of resonances \cite{Benettin}. Second, if action variables can leave
the initial cylinder in an exponentially long time, they must enter
some resonant block associated to a lower dimensional lattice.

Obviously, the condition of nonoverlapping of resonances is equivalent to the following form:
$$|k\cdot\widetilde{\omega}(x)+2\pi l_0|> \alpha_{r} \quad\text{for~ all}\,~x\in\mathcal {V}_{\delta_r}(\mathcal {B}_{\mathcal {M},\,\delta_r}^l),~ k\in\mathbb{Z}_{K}^n\backslash\mathcal {M}~\text{and} ~l_0\in \mathbb{Z}\,.$$
In order to satisfy the above condition, we can make choice
\begin{align*}
\alpha_r&=\Big(\frac{2}{\mu}\Big)^rr!K^{\frac{r(r-1)}{2}}\alpha_1\,, \quad r\geq 2\\
\delta_r&=\frac{\alpha_r}{3K\beta}\,.
\end{align*}

Indeed, one can remark that for any $x\in\mathcal {V}_{\delta_r}(\mathcal {B}_{\mathcal {M},\,\delta_r}^l)$, there exists $x_*\in \mathcal {B}_{\mathcal {M}}^l$ such that $|x-x_*|_2\leq \delta_r+\frac{4}{\mu}\delta_r+\frac{2rK^{r-1}}{\beta\mu}\alpha_r$ by the definition of the extended block and Proposition 3.1 (iv). Therefore, for any $x\in\mathcal {V}_{\delta_r}(\mathcal {B}_{\mathcal {M},\,\delta_r}^l),~k\in\mathbb{Z}_{K}^n\backslash\mathcal {M}$\,, and $l_0\in\mathbb{Z}$, we have
\begin{align*}
|k\cdot\widetilde{\omega}(x)+2\pi l_0|&\geq |k\cdot\widetilde{\omega}(x_*)+2\pi l_0|-|k\cdot(\widetilde{\omega}(x)-\widetilde{\omega}(x_*))|\\
&>\alpha_{r+1}-K\beta|x-x_*|_2\\
&\geq \alpha_{r+1}-K\beta\Big(\delta_r+\frac{4}{\mu}\delta_r+\frac{2rK^{r-1}}{\beta\mu}\alpha_r\Big)\,.
\end{align*}

Thus the nonoverlapping condition of resonances is satisfied only if $\alpha_{r+1}-K\beta\big(\delta_r+\frac{4}{\mu}\delta_r+\frac{2rK^{r-1}}{\beta\mu}\alpha_r\big)\geq\alpha_r$ through the above choices. By means of the nonoverlapping condition of resonances, we can prove that
\begin{equation}\label{eq:25}
\Big|1-e^{ik\cdot\widetilde{\omega}(x)}\Big|\geq\frac{2\alpha_r}{\pi}\quad \text{for~all}\,~x\in\mathcal {V}_{\delta_r}(\mathcal {B}_{\mathcal {M},\,\delta_r}^l),~k\in\mathbb{Z}_{K}^n\backslash\mathcal {M}\,.
\end{equation}
For this purpose, we can choose $l=l(\widetilde{\omega},k)\in\mathbb{Z}$ such that $\Big|\frac{k\cdot\widetilde{\omega}(x)+2\pi l}{2}\Big|\in \big[0,\frac{\pi}{2}\big]$. Thus,
$$\Big|1-e^{ik\cdot\widetilde{\omega}(x)}\Big|=2\sin\bigg|\frac{k\cdot\widetilde{\omega}(x)+2\pi l}{2}\bigg|\geq \frac{4}{\pi}\bigg|\frac{k\cdot\widetilde{\omega}(x)+2\pi l}{2}\bigg|\geq\frac{2}{\pi}\alpha_r\,.$$

When choosing $\rho_1^{(r)}=\min\big(\frac{\alpha_r}{4K\beta},1\big)<\delta_r$ and $\rho_2^{(r)}=\widetilde{\sigma}_2$, Normal Form Lemma can be applied in the domain $\mathcal {V}_{\rho_1^{(r)}/2}(\mathcal {B}_{\mathcal {M},\,\delta_r}^l)\times\mathcal {V}_{\rho_2^{(r)}/2}(T^n)$ with $\alpha=\frac{2}{\pi}\alpha_r$, if the following conditions are satisfied:
\begin{align}
||Df||_{G^{(r)},\,\frac{\rho^{(r)}}{2},\,c^{(r)}}&\leq \frac{2}{\pi}\cdot\frac{\alpha_r\rho_1^{(r)}}{\widetilde{C}AK\rho_2^{(r)}}\,,\label{eq:16}\\
\rho_1^{(r)}&\leq \widetilde{\sigma}_1\,,\label{eq:17}
\end{align}
where $G^{(r)}=\mathcal {B}_{\mathcal {M},\,\delta_r}^l,~c^{(r)}=\frac{\rho_1^{(r)}}{\rho_2^{(r)}}$\,.

On the other hand, by Lemma 2.3 we have $||Df||_{G^{(r)},\,\frac{\rho^{(r)}}{2},\,c^{(r)}}\leq\frac{2}{\rho_2^{(r)}}||f||_{\mathcal {G}_{\gamma},\,\widetilde{\sigma}}=\frac{2}{\rho_2^{(r)}}\tilde{\epsilon}$. Therefore, in order to satisfy (\ref{eq:16}), we require
\begin{equation}\label{eq:18}
\tilde{\epsilon}\leq\frac{\alpha_r\rho_1^{(r)}}{\widetilde{C}AK\pi}\,.
\end{equation}

According to Normal Form Lemma, there exists a symplectic transformation $\Phi :\mathcal
{D}_{\rho^{(r)}/4}(G^{(r)})\rightarrow \mathcal {D}_{\rho^{(r)}/2}(G^{(r)})$ such that
the conjugate symplectic map $T_{\gamma}'=\Phi^{-1}\circ
T_{\gamma}\circ\Phi$ is generated by the analytic function
$F'=F_0+Z+R$ with $Z\in \mathcal {R}(\mathcal {M},K)$ and
$$||DR||_{G^{(r)},\,\frac{\rho^{(r)}}{4},\,c^{(r)}}\leq 3e^{-\frac{K\rho_2^{(r)}}{24}}\frac{2}{\rho_2^{(r)}}\cdot\tilde{\epsilon}\,.$$

\textbf{Claim 1.} Under the above conditions, denote by $\tau_{+}$ and $\tau_{-}$ the possibly times of escape of $x(t)$ from $\mathcal {Z}_{\mathcal {M}}^l$ at positive and negative times respectively. Then for any $x_0\in\mathcal {B}_{\mathcal {M}}^l$\,, one has $x(t)\in\mathcal {C}_{\mathcal {M},\,\delta_r}^l(x_0)$ if $|t|\leq\min(T^{(r)},\,\tau_{+},\,-\tau_{-})$ with $T^{(r)}=\frac{\rho_1^{(r)}\rho_2^{(r)}e^{\frac{K\rho_2^{(r)}}{24}}}{2^7\tilde{\epsilon}}$\,.

Consider the action variables of $T_{\gamma}'$, it is given by
$$\hat{a}=a-\Big(\frac{\partial Z(\hat{a},\varphi)}{\partial \varphi}+\frac{\partial R(\hat{a},\varphi)}{\partial \varphi}\Big)\,,$$
and consider a auxiliary system
$$\tilde{a}=a-\frac{\partial Z(\hat{a},\varphi)}{\partial \varphi}\,.$$
After iterating $t$ times, we have $\tilde{a}(t)\in\Pi_{\mathcal {M}}(a(0))$, and
$$|a(t)-\tilde{a}(t)|_2\leq |t|\cdot\Big|\Big|\frac{\partial R(\hat{a},\varphi)}{\partial \varphi}\Big|\Big|_{G^{(r)},\,\frac{\rho^{(r)}}{4},\,1}\leq T^{(r)}||DR||_{G^{(r)},\,\frac{\rho^{(r)}}{4},\,c^{(r)}}\leq\frac{\rho_1^{(r)}}{2^4}\,.$$
Thus,
\begin{align*}
\text{dist}(x(t),\Pi_{\mathcal {M}}(x(0)))&\leq |x(t)-a(t)|_2+|a(t)-\tilde{a}(t)|_2+\text{dist}(\tilde{a}(t),\Pi_{\mathcal {M}}(x(0)))\\
&\leq \frac{\rho_1^{(r)}}{2^8}+\frac{\rho_1^{(r)}}{2^4}+\frac{\rho_1^{(r)}}{2^8}\\
&\leq\frac{\rho_1^{(r)}}{4}\,.
\end{align*}
That is $x(t)\in\mathcal {C}_{\mathcal {M},\,\delta_r}^l(x_0)$ if $|t|\leq\min(T^{(r)},\,\tau_{+},\,-\tau_{-})$.

\textbf{Claim 2.} Assume $\tau_{+},-\tau_{-}<T^{(r)}$ and there exists $t_*$ satisfying $x(t_*)\in\mathcal {C}_{\mathcal {M},\,\delta_r}^l(x_0)$, but $x(t_*+1)\notin\mathcal {C}_{\mathcal {M},\,\delta_r}^l(x_0)$, then $x(t_*+1)\in \mathcal {B}_{\mathcal {M}'}^{l'}$ with dim$(\mathcal {M}')<r$ and some $l'$.

Note that
\begin{align*}
|x(t_*+1)-x(t_*)|_2&\leq|x(t_*+1)-a(t_*+1)|_2+\Big|\Big|\frac{\partial Z(\hat{a},\varphi)}{\partial \varphi}\Big|\Big|_{G^{(r)},\,\frac{\rho^{(r)}}{4},\,1}\\
&\quad+\Big|\Big|\frac{\partial R(\hat{a},\varphi)}{\partial \varphi}\Big|\Big|_{G^{(r)},\,\frac{\rho^{(r)}}{4},\,1}+|a(t_*)-x(t_*)|_2\\
&\leq \frac{\rho_1^{(r)}}{2^8}+2||Df||_{G^{(r)},\,\frac{\rho^{(r)}}{2},\,c^{(r)}}+3||Df||_{G^{(r)},\,\frac{\rho^{(r)}}{2},\,c^{(r)}}+\frac{\rho_1^{(r)}}{2^8}\\
&\leq \frac{\rho_1^{(r)}}{2^7}+5\cdot\frac{2}{\rho_2^{(r)}}\tilde{\epsilon}
\end{align*}
where we have used Normal Form Lemma. Because of $\widetilde{C}\geq 51,~A\geq 1,~\alpha_r<1$ and (\ref{eq:18}), we have
$$\frac{2}{\rho_2^{(r)}}\tilde{\epsilon}\leq \frac{\rho_1^{(r)}}{6*51}\cdot\frac{1}{K\rho_2^{(r)}}\leq \frac{\rho_1^{(r)}}{6*51}\,,$$
when $K\geq\frac{1}{\rho_2^{(r)}}=\frac{1}{\widetilde{\sigma_2}}$\,. Thus, $|x(t_*+1)-x(t_*)|_2\leq\frac{\rho_1^{(r)}}{4}<\frac{\delta_r}{4}$\,.
That means $x(t_*+1)\in\mathcal {V}_{\delta_r}(\mathcal {B}_{\mathcal {M},\,\delta_r}^l)$. At this time there are only two possible cases: $x(t_*+1)\in\mathcal {Z}_{\mathcal {M}}^{l'}~\text{for some}~l'\in\mathbb{Z}^r$ or $x(t_*+1)\notin \mathcal {Z}_r^*$ by the nonoverlapping condition of resonances. However, for any $l'\neq l$, $x(t_*+1)\notin\mathcal {Z}_{\mathcal {M}}^{l'}$. In fact, because $x(t_*)\in\mathcal {Z}_{\mathcal {M}}^{l}$, we have $\big|k^{(j)}\cdot\widetilde{\omega}(x(t_*))+2\pi l_j\big|\leq \alpha_r,~j=1,\ldots,r$. Thus, for any $l_j'\in \mathbb{Z}$ and $l_j'\neq l_j$,
\begin{align*}
\big|k^{(j)}\cdot\widetilde{\omega}(x(t_*+1))+2\pi l_j'\big|&=\big|k^{(j)}\cdot\big(\widetilde{\omega}(x(t_*+1))-\widetilde{\omega}(x(t_*))\big)+k^{(j)}\cdot\widetilde{\omega}(x(t_*))+2\pi l_j\\
&\quad+2\pi l_j'-2\pi l_j\big|\\
&\geq 2\pi-\alpha_r-K\beta\cdot\frac{\delta_r}{4}\\
&> 1.
\end{align*}

Therefore $x(t_*+1)\notin \mathcal {Z}_r^*$. Because of the property (iii) in Proposition 3.1, one has $x(t_*+1)\in\mathcal {B}_{\mathcal {M}'}$ with dim$(\mathcal {M}')<r$. It is shown that $x(t)$ enter some resonant block associated to lower dimensional lattice, after going out of the original cylinder only through its base during $|t|<T$.

The stability estimates now apply to all blocks simultaneously, if
\begin{gather}
\tilde{\epsilon}=\frac{M}{\beta}\cdot\epsilon\leq\frac{\alpha_1\rho_1^{(1)}}{\widetilde{C}AK\pi}\,,\label{eq:19}\\
\rho_1^{(n)}\leq \widetilde{\sigma}_1\,,\label{eq:20}\\
\beta\leq \min_{1\leq i\leq n}\{\alpha_i\}\,,\label{eq:21}\\
K\geq \frac{1}{\widetilde{\sigma}_2}=\frac{1}{\sigma_2}\,.\label{eq:22}
\end{gather}
The (\ref{eq:19}) -- (\ref{eq:20}) are satisfied if we require
\begin{align}
&\epsilon\leq\frac{\alpha_n^2}{4\pi(21n+30)\big[(\frac{2M}{m})^n\cdot n!\big]^2K^{n^2-n+2}\cdot M}\,,\label{eq:23}\\
&\alpha_n\leq 4KM\sigma_1\,.\label{eq:24}
\end{align}

On the other hand, in order to keep as small as possible for the diameter of the cylinders, at least the order of $K^{-1}$, it is convenient to put $\alpha_n= 4MK^{-n}\sigma_1$. Due to $\sigma_1<\frac{1}{4M}$, so $\alpha_n<1$. Finally, let $\beta=\min\limits_{1\leq i\leq n}\{\alpha_i\}$ and $K=\epsilon^{-\frac{1}{b}}$ with $b=2(n^2+n+2)$. Then, $K\geq\frac{1}{\sigma_2}$ when $\epsilon\leq\sigma_2^b$.

From the above discussions, it is clear that for any initial value $x_0$, the adapted normal form can be constructed in the corresponding extended resonant block. The normal form provides the confinement of the action $x(t)$ in $\Pi_{\mathcal {M}}(x_0)$. Moreover, $x(t)$ must enter one of the other resonant block with lower dimensional multiplicity after out of the previous one, arrives in the nonresonant block at worst, where it stops. Thus, the stability radius $\widetilde{\triangle}$ satisfies
$$\widetilde{\triangle}\leq (n-1)(\frac{4M}{m}\delta_n+\frac{2nK^{n-1}M}{\beta m}\alpha_n)+\widetilde{\triangle}_0\,,$$
where $\widetilde{\triangle}_0$ denotes the deviation of the action variables when them get into the nonresonant block and $|t|\leq T^{(1)}$. This number can be estimated directly as follows.

According to Proposition 3.1 (i), for any $x\in\mathcal {B}_0$, $\big|k\cdot\widetilde{\omega}(x)+2\pi l\big|>\alpha_1$ for all $k\in\mathbb{Z}_{K}^n\backslash\{0\}$ and $l_0\in \mathbb{Z}$. Then we can prove
$$\Big|1-e^{ik\cdot\widetilde{\omega}(x)}\Big|\geq\frac{2}{\pi}\alpha_1\quad \text{for~all}\,~x\in\mathcal {B}_0,~k\in\mathbb{Z}_{K}^n\backslash\{0\}\,.$$
by the same technique that be used in the proof of (\ref{eq:25}). Furthermore, applying Lemma 2.1 with $G=\mathcal {B}_0,~\rho_1=\rho_1^{(1)},~\alpha=\frac{2}{\pi}\alpha_1$, and $\mathcal {M}=\{0\}$, we have
$$\Big|1-e^{ik\cdot\widetilde{\omega}(x)}\Big|\geq\frac{\alpha_1}{\pi}\quad \text{for~all}\,~x\in\mathcal {V}_{\rho_1^{(1)}}(\mathcal {B}_0),~k\in\mathbb{Z}_{K}^n\backslash\{0\}\,.$$

Now, we can apply Normal Form Lemma with $G=\mathcal {B}_0,~\rho=\frac{\rho^{(1)}}{2},~\alpha=\frac{\alpha_1}{\pi}$, and $\mathcal {M}=\{0\}$ under the condition (\ref{eq:23}). Then there exists a real analytic canonical transformation $\Phi :\mathcal
{D}_{\rho^{(1)}/4}(\mathcal {B}_0)\rightarrow \mathcal {D}_{\rho^{(1)}/2}(\mathcal {B}_0)$ such that $T_{\gamma}'=\Phi^{-1}\circ
T_{\gamma}\circ\Phi$ is generated by the analytic function
$F'(\hat{a},\varphi)=F_0(\hat{a})+Z(\hat{a},\varphi)+R(\hat{a},\varphi)$ with $Z\in \mathcal {R}(\mathcal {M},K)$. Due to $\mathcal {M}=\{0\}$, $Z$ only depends on the action variables, and the system becames
\begin{equation*}
\left\{ \begin{aligned}
        \hat{a}&=a-\frac{\partial R(\hat{a},\varphi)}{\partial \varphi} \\
        \hat{\varphi} &=\varphi+\frac{\partial F_0(\hat{a})}{\partial \hat{a}}+\frac{\partial Z(\hat{a})}{\partial \hat{a}}+\frac{\partial R(\hat{a},\varphi)}{\partial \hat{a}}\,.
        \end{aligned} \right.
\end{equation*}
After iterating $t$ times with $|t|\leq T^{(1)}$,
\begin{align*}
|a(t)-a(0)|_2&\leq |t|\cdot\Big|\Big|\frac{\partial R}{\partial \varphi}\Big|\Big|_{\mathcal {B}_0,\,\frac{\rho^{(1)}}{4},\,1}\\
&\leq T^{(1)}\cdot 3e^{-\frac{K\rho_2^{(1)}}{24}}||Df||_{\mathcal {B}_0,\,\frac{\rho^{(1)}}{2},\,c^{(1)}}\\
&\leq \frac{\rho_1^{(1)}}{2^4}\,.
\end{align*}
Thus, $|x(t)-x(0)|_2\leq |x(t)-a(t)|_2+|a(t)-a(0)|_2+|a(0)-x(0)|_2\leq \frac{\rho_1^{(1)}}{4}\leq \delta_n$\,. That means $\widetilde{\triangle}_0\leq\delta_n$ and
$$\widetilde{\triangle}\leq n(\frac{4M}{m}\delta_n+\frac{2nK^{n-1}M}{\beta m}\alpha_n)\,.$$
Therefore, we have
$$|x(t)-x(0)|_2\leq n(\frac{4M}{m}\delta_n+\frac{2nK^{n-1}M}{\beta m}\alpha_n)\quad \text{for~all}\,~x(0)\in\mathcal {G}_{\gamma}~\text{and}\,~|t|\leq T^{(1)}\,.$$
Note that $I=\gamma x$, so
$$
|I(t)-I(0)|_2\leq \gamma\widetilde{\triangle}\leq n\Big(\frac{16M}{3m}+\frac{8nM}{m}\Big)\sigma_1\cdot\epsilon^{\frac{1}{b}}:=c_0\epsilon^{\frac{1}{b}}:=\triangle\,.
$$
where $c_0=\frac{8nM}{3m}(3n+2)\sigma_1$.

Combining the above discussion, we can choose the stability time $T=\min\limits_{1\leq i\leq n}\{T^{(i)}\}=T^{(1)}$. Precisely, $$T=\frac{\rho_1^{(1)}\rho_2^{(1)}e^{\frac{K\rho_2^{(1)}}{24}}}{2^7\tilde{\epsilon}}
=\frac{\alpha_1\sigma_2}{2^9\epsilon KM}\cdot e^{\frac{\sigma_2}{24}\epsilon^{-\frac{1}{b}}}
=T_0\epsilon^{-\frac{3}{4}}e^{\frac{\sigma_2}{24}\epsilon^{-\frac{1}{b}}}\,,$$
with $T_0=\frac{(\frac{m}{2M})^n\sigma_1\sigma_2}{2^7n!}$\,. And the perturbation $\epsilon\leq\min(\epsilon_0,\,\sigma_2^b)$, where $$\epsilon_0=\frac{M^2\cdot\sigma_1^4}{\pi^2(21n+30)^2\big[(\frac{2M}{m})^n\cdot n!\big]^4}\,.$$

Finally, in order to prevent action variables from going out of $\mathcal {G}$, we restrict $I_0\in \mathcal {G}-\triangle$. 

\section{Application}

An application of the above theorem gives the exponential stability of a nearly integrable symplectic map with a small twist, which often comes from numerical discretization of Hamiltonian systems. Consider a one-parameter family of symplectic map $\mathcal {C}_s$ with the parameter $s$ satisfying $0<s\leq 1$ and the analytic generating function $H_s=sH_0(\hat{I})+sh(\hat{I},\theta)$. The small twist map is given by
\begin{equation*}
 \begin{aligned}
         \hat{I} &= I-s\frac{\partial h(\hat{I},\theta)}{\partial \theta}~, \\
         \hat{\theta} &=\theta+s\frac{\partial H_0(\hat{I})}{\partial \hat{I}}+s\frac{\partial h(\hat{I},\theta)}{\partial \hat{I}}~.
                          \end{aligned}
                          \end{equation*}
The result can ba stated as follows.\\[-5pt]

\textbf{Theorem 5.1} Let
$H_s$ be analytic in
$\mathscr{D}_{\sigma}(\mathcal {G})$, where $\mathcal {G}$ is an
open bounded domain of $\mathbb{R}^n$, $\sigma=(\sigma_1,\sigma_2)$ is positive, and
$\omega(\hat{I})=\partial H_0(\hat{I})$ satisfying (\ref{eq:30}). Consider the above small twist symplectic map $\mathcal {C}_s$ defined on $\mathcal {G}\times T^n$. If
\begin{equation}\label{eq:29}
\epsilon=||h||_{\mathcal {G},\,\sigma}\leq \min(\epsilon_0,\,\sigma_2^b)\,,
\end{equation}
where
$$\epsilon_0=\frac{M^2\cdot\sigma_1^4}{\pi^2(21n+30)^2\big[(\frac{2M}{m})^n\cdot n!\big]^4}\,,$$
$b=2(n^2+n+2)$ and $\sigma_1< \frac{1}{4M}$\,.
Then for any $0<s\leq1$, the symplectic map $\mathcal {C}_s$ satisfies
$$|I(t)-I_0|_2\leq \triangle \quad \mbox{for}\quad |t|\leq T \quad \mbox{and}\quad I_0\in \mathcal {G}-\triangle\,,$$
where $(I(t),\theta(t))=\mathcal {C}^t(I_0,\theta_0)$ ($t$ is viewed
as iterative times),

\begin{equation*}
 \begin{aligned}
         \triangle &= c_0\epsilon^{\frac{1}{b}}\quad&\mbox{with}\quad &c_0=\frac{8nM}{3m}(3n+2)\sigma_1\,,\\
         T &=T_0\epsilon^{-\frac{3}{4}}e^{c_1\epsilon^{-\frac{1}{b}}}\quad&\mbox{with}\quad
         &T_0=\frac{(\frac{m}{2M})^n\sigma_1\sigma_2}{2^7n!}\quad\mbox{and}\quad
         c_1=\frac{\sigma_2}{24}\,.
                          \end{aligned}
                          \end{equation*}\\[-5pt]

Here we only outline the proof of the above theorem since it is almost the same as the proof of Theorem 1. First, in this case, the small denominators became $1-e^{ik\cdot s\widetilde{\omega}(x)}$. Therefore, in order to construct the resonant normal form with respect to the $K$-lattice $\mathcal {M}$, we restrict to the subset
$$G_s=\big\{x\in\mathcal {G}_{\gamma}:~|1-e^{ik\cdot s\widetilde{\omega}(x)}|\geq s\alpha\,,\quad \forall~\,k\in\mathbb{Z}_K^n\backslash\mathcal {M}\big\}\,.$$

The corresponding Normal Form Lemma holds if the quantity are $sf,\,sM,\,sm,\,s\alpha$ and $s\beta$ instead of $f,\,M,\,m,\,\alpha$ and $\beta$ respectively. That means if $s||Df||_{G,\,\rho,\,c}\leq \frac{s\alpha\,\rho_1}{\widetilde{C}AK\rho_2}\,,$ then there exists a real analytic
canonical transformation $\Phi :\mathcal {D}_{\frac{\rho}{2}}(G)\rightarrow \mathcal {D}_{\rho}(G)$ such that
the conjugate symplectic map $T_{\gamma}'=\Phi^{-1}\circ T_{\gamma}\circ\Phi:\mathcal {D}_{\frac{\rho}{2}}(G)\rightarrow
\mathcal {D}_{\rho}(G)$ is generated by the analytic function $F'=sF_0+sZ+sR$ with $Z\in \mathcal {R}(\mathcal {M},K)$ and $R$ exponential small. Moreover, the same estimates about $Z$ and $R$ hold like in Lemma 2.2.

On the other hand, the geometric construction also needs some modifications. Precisely, the parameter $s$ enters the definitions of resonant manifold, zone, and block. $\widetilde{\omega}(x)$ and $\alpha_r$ should be replaced by $s\widetilde{\omega}(x)$ and $s\alpha_r$ respectively. Thus we obtain the corresponding resonant manifold, zone, and block. We remark that the arguments of the proposition 3.1 still hold in this case.

Finally, we make the same choices about the parameters $\alpha_r,\delta_r~(r=1,\ldots,n)$ and $K$ as before, and the desired estimates can be derived. The details are omitted here.

In the following contents, we will discuss the exponential stability of symplectic integrator which applied to the integrable Hamiltonian system by the Theorem 5.1. Integrable Hamiltonian system is a very important class of dynamic system. In general it possesses many enough first integrals. Therefore, it exhibits regular dynamic behavior which corresponds to periodic and quasi-periodic motions in phase space through action-angle variables. However, in many cases the action-angle variables may not be known explicitly. Then it is difficult to compute solutions of the given integrable Hamiltonian system, and numerical integration is necessary.

After the pioneering work of Channel (1983), Feng Kang (1985, 1986) and Ruth (1983), the symplectic integrator has become a widely interested subject on the problem of numerically solving Hamiltonian systems. Extensive computer experimentation, by some typical models of Hamiltonian systems,
has shown the overwhelming superiority of symplectic algorithms over the
conventional non-symplectic ones, especially in simulating the global and
structural dynamic behavior of the systems (e.g. see \cite{Haire} and \cite{FQ}). The symplectic algorithm, which applied to integrable Hamiltonian system, may be characterized as a perturbation of the phase flow of the integrable system. Here the smallness of the perturbation is described by the
time-step size of the algorithm which also enters into the frequency map of the integrable system. Therefore numerical stability problem arises.

There has been recently some nice work about the numerical analysis of symplectic algorithms for Hamiltonian systems, for example,
by Benettin \& Giorgilli \cite{BG}, Hairer \& Lubich \cite{HL}, Shang \cite{Shang1999}, and Stoffer \cite{St}. Stoffer proved the numerical solutions is integrable up to a remainder which is exponential small with respect to the step-size when a symplectic integrator is applied to a integrable system. However, the result requires that the initial frequency satisfying the strong non-resonance condition. In addition, for non-resonance time step-size, Shang obtained the existence result of numerical invariant tori of symplectic algorithms.

We consider a integrable Hamiltonian system (usually not given in action-angle variables)
\begin{equation}\label{eq:26}
\dot{p}=-\frac{\partial H(p,q)}{\partial q}, \quad\dot{q}=\frac{\partial H(p,q)}{\partial p}
\end{equation}
and apply to it a symplectic algorithm $G_H^h$ of order $r$ with step size $h$. For an overview on the symplectic integrators see the book of Haire, Lubich and Wanner \cite{Haire}. Because of Arnold-Liouville theorem, there exists a symplectic transformation $\Psi :(I,\theta)\rightarrow (p,q)$ such that the new Hamiltonian $\mathcal {H}(I)=H \circ\Psi(I,\theta)$, that only depends on the action variables, $(I,\theta)\in \mathcal {G}\times T^n$. Here we assume $H$ and $\mathcal {H}$ is analytic in $\mathcal {D}_{\rho}(\mathcal {G})$ and $\mathcal {V}_{\rho_1}(\mathcal {G})$ respectively. And $\omega(I)=\frac{\partial \mathcal {H}(I)}{\partial I}$ satisfies the condition (\ref{eq:30}). In the action-angle variables, the equation (\ref{eq:26}) takes the simple form
\begin{equation}\label{eq:27}
\dot{I}=0, \quad\dot{\theta}=\omega(I)=\frac{\partial \mathcal {H}(I)}{\partial I}.
\end{equation}
The symplectic integrator $G_H^h$ becames $\widetilde{G}_H^h=\Psi^{-1}\circ G_H^h\circ\Psi$.\\[-5pt]

\textbf{Lemma 5.1} \cite{Shang1999} There exists a function $f^h$ which depends on the time step $h$
such that it is well-defined and real analytic in the domain $\mathcal {D}_{\frac{\rho}{4}}(\mathcal {G})$ for $h\in[0,\delta]$ with $\delta$ being a sufficiently small positive number so that $\widetilde{G}_H^h:(I,\theta)\rightarrow (I',\theta')$ can be expressed by $f^h$ as follows:
\begin{equation}\label{eq:28}
I'=I-h^{r+1}\frac{\partial f^h(I',\theta)}{\partial \theta}, \quad\theta'=\theta+h\omega(I')+h^{r+1}\frac{\partial f^h(I',\theta)}{\partial I'}.
\end{equation}
Moreover, there exists $L$ independent on $h$ such that $||f^h||\leq L$.\\[-5pt]

Now the Theorem 5.1 can be applied to $\mathcal {C}_h=\widetilde{G}_H^h$ if $h^rf^h$ satisfies the estimate (\ref{eq:29}) with $\sigma=\frac{\rho}{4}$. Therefore we have the following result.\\[-5pt]

\textbf{Theorem 5.2} Under the above assumption on $\mathcal {H}$, apply a symplectic method $G_H^h$ of order $r$ to the equation (\ref{eq:26}). Let $\widetilde{G}_H^h$ generate an orbit $(I_1,\theta_1),(I_2,\theta_2),\ldots$ with any initial value $(I_0,\theta_0)$ in action-angle variables. Then there are positive constants $h_0,c_0,c_1,T_0$ such that for all $h\leq h_0$, the following estimates hold
$$|I_m-I_0|_2\leq c_0h^{\frac{r}{b}}\,,$$
for all $m$ with $mh\leq T_0h^{-\frac{3}{4}r}e^{c_1h^{-\frac{r}{b}}}$ and $b=2(n^2+n+2)$.


\begin{thebibliography}{1}

\bibitem{Arnold}Arnold, V.I. (1963). Proof of A. N. Kolmogorov's theorem on the
preservation of quasi-periodic motions under small perturbations of
the Hamiltonian. Russ. Math. Surv., Vol.18(5), pp. 9-36.

\bibitem{Arnold-2}Arnold V.I. (1989). Mathematical Methods of Classical Mechanics. 2nd
ed. Springer- Verlag, New York.

\bibitem{Benettin}Benettin G., Galgani L. and Giorgilli A. (1985). A proof of Nekhoroshev¡¯s theorem
for the stability times in nearly integrable Hamiltonian systems. Cel. Mech., Vol. 37, pp. 1-25.

\bibitem{BG}Benettin G., Giorgilli A. (1994). On the Hamiltonian interpolation of near to the identity
symplectic mappings with application to symplectic integration algorithms. J. Statist.
Phys., Vol. 74, pp. 1117-1143.

\bibitem{AB}Bounemoura A. and Marco J.-P. (2011). Improved exponential stability for
near-integrable quasi-convex Hamiltonians. Nonlinearity, Vol. 24(1), pp. 97-112.

\bibitem{Channell83}Channell P.J. (1983). Symplectic integration algorithms. Los Alamos National Laboratory
Report AT-6:ATN-83-9.

\bibitem{Channell90}Channell P.J., Scovel C. (1990). Symplectic integration of Hamiltonian systems. Nonlinearity
3, 231-259

\bibitem{Cheng}Cheng C. (2006). Hamiltonian Systems: Stable or Unstable? Milan j. math., Vol. 74, pp. 295¨C312.

\bibitem{Delshams-Gutierrez}Delshams A. and Gutierrez P. (1996).
Effective stability and KAM theory. J. Diff. Eq., Vol. 128,
pp. 415-490.

\bibitem{Ding}Ding Z. and Shang Z. (2018) KAM invariant tori of symplectic integrators for R\"{u}ssmann's non-degenerate Hamiltonian systems. Accepted by SCIENCE CHINA Mathematics.

\bibitem{FQ}Feng K., Qin M. (2003). Symplectic Geometric Algorithms For Hamiltonian System, Zhejiang Science \& Technology Press, Hangzhou

\bibitem{Guzzo}Guzzo M. (2004). A direct proof of the Nekhoroshev theorem for nearly
integrable sysmplectic maps. Ann. Henri Poincar\'{e}, Vol. 5, pp.
1013-1039.

\bibitem{HL}Hairer E., Lubich C. (1997). The life-span of backward error analysis for numerical
integrators. Numer. Math., Vol. 76(4), pp. 441-462.

\bibitem{Haire}Hairer E., Lubich C. and Wanner G. (2006). Geometric Numerical Integration: Structure-Preserving Algorithms for Ordinary Differential
Equations, 2nd ed. Springer- Verlag, Berlin.

\bibitem{Kuksin-Poschel}Kuksin S. B. and P\"{o}schel J. (1994). On the inclusion of analytic
symplectic maps in analytic Hamiltonian flows and its applications.
Nonlinear Differential Equations Appl., Vol. 12, pp. 96-116.

\bibitem{Lochak-1}Lochak P., Neishtadt A. I. (1992). Estimates of stability time for nearly
integrable systems with a quasiconvex Hamiltonian. Chaos, Vol. 2,
pp. 492-499.

\bibitem{Lochak-2}Lochak P. (1992). Canonical perturbation theory via simultaneous approximation. Russian
Math. Surveys, Vol. 47, pp. 57-133.

\bibitem{Nekhoroshev}Nekhoroshev N. N. (1977). An exponential estimate of the time of stability of nearly integrable
Hamiltonian systems I. Uspekhi Mat. Nauk, Vol. 32, pp. 5-66.

\bibitem{Poschel}P\"{o}schel J. (1993). Nekhoroshev estimates for quasi-convex Hamiltonian
systems. Math. Z., Vol. 213, pp. 187-216.

\bibitem{Shang1999}Shang Z. J. (1999). On the KAM theorem of symplectic algorithms for Hamiltonian
systems. Numer. Math., Vol. 83, pp. 477-496.

\bibitem{Shang2000}Shang Z. J. (2000). A note on the KAM theorem for symplectic
mappings. J. Dynam. Differential Equations, Vol. 12(2), pp. 357-383.

\bibitem{St}Stoffer D. (1998). On the qualitative behaviour of symplectic integrators. Part II. Integrable systems. J. Math. Anal. Appl., Vol. 217, pp. 501-520.

\bibitem{XieBo}Xie B. and Shang Z. Effctive stability analysis for nearly integrable
symplectic maps with applications to symplectic integrators. Unpublished




















\end{thebibliography}
\end{document}